\makeatletter 

\expandafter\ifx\csname amscd.sty\endcsname\relax
\expandafter\def\csname amscd.sty\endcsname{}
\else\endinput\fi
\def\filename{amscd.sty}
\def\fileversion{1.1} \def\filedate{21-JUN-1991}
\immediate\write16{%
AMS-Latex option `\filename' (\fileversion, \filedate)}
\def\Invalid@@{Invalid use of \string}
\def\Let@{\let\\\math@cr}
\def\RIfM@{\relax\protect\ifmmode}
\@ifundefined{math@cr}
  {\def\math@cr{{\ifnum0=`}\fi
   \new@ifstar{\global\@eqpen\@M\math@cr@}%
          {\global\@eqpen\interdisplaylinepenalty \math@cr@}}}
  {}
\def\math@cr@{\new@ifnextchar[\math@cr@@{\math@cr@@[\z@]}}
\def\math@cr@@[#1]{\ifnum0=`{\fi}\math@cr@@@
  \noalign{\vskip#1\relax}}
\def\restore@math@cr{\def\math@cr@@@{\cr}}
\restore@math@cr
\def\new@ifnextchar#1#2#3{%
  \let\@tempe #1\def\@tempa{#2}\def\@tempb{#3}\futurelet
    \@tempc\new@ifnch}
\def\new@ifnch{\ifx\@tempc \@tempe \let\@tempd\@tempa
             \else\let\@tempd\@tempb\fi\@tempd}
\def\new@ifstar#1#2{\new@ifnextchar *{\def\@tempa*{#1}\@tempa}{#2}}
\def\DN@{\def\next@}
\def\FN@{\futurelet\next}
\def\setboxz@h{\setbox\z@\hbox}
\def\wdz@{\wd\z@}
\def\setbox@ne{\setbox\@ne}
\def\wd@ne{\wd\@ne}
\def\rightarrowfill@#1{\m@th\setboxz@h{$#1-$}\ht\z@\z@
  $#1\copy\z@\mkern-6mu\cleaders
  \hbox{$#1\mkern-2mu\box\z@\mkern-2mu$}\hfill
  \mkern-6mu\mathord\rightarrow$}
\def\leftarrowfill@#1{\m@th\setboxz@h{$#1-$}\ht\z@\z@
  $#1\mathord\leftarrow\mkern-6mu\cleaders
  \hbox{$#1\mkern-2mu\copy\z@\mkern-2mu$}\hfill
  \mkern-6mu\box\z@$}
\def\leftrightarrowfill@#1{\m@th\setboxz@h{$#1-$}\ht\z@\z@
  $#1\mathord\leftarrow\mkern-6mu\cleaders
  \hbox{$#1\mkern-2mu\box\z@\mkern-2mu$}\hfill
  \mkern-6mu\mathord\rightarrow$}
\long\def\@leftmark#1#2{#1}
\long\def\@rightmark#1#2{#2}
\long\def\@ifempty#1{%
 \expandafter\ifx\@car#1@\@nil @\@empty
  \expandafter\@leftmark\else\expandafter\@rightmark\fi}
\long\def\@ifnotempty#1{\@ifempty{#1}{}}
\def\atdef@#1{\expandafter\def\csname\space @\string#1\endcsname}
\begingroup \catcode`\@=\active
\xdef @{\expandafter\noexpand\csname FN\string @\endcsname
  \expandafter\noexpand\csname at\string @\endcsname}
\endgroup
\def\at@{\let\next@\at@@
 \ifcat\noexpand\next a\else
 \ifcat\noexpand\next0\else
 \ifcat\noexpand\next\relax\else
 \let\next@\at@@@\fi\fi\fi\next@}
\def\at@@#1{\expandafter
  \ifx\csname\space @\string#1\endcsname\relax
    \DN@{\at@@@#1}%
  \else
    \DN@{\csname\space @\string#1\endcsname}%
  \fi\next@}%
\def\at@@@{\err@{\Invalid@@ @}{\the\athelp@}\char64\relax}
\@ifundefined{athelp@}{\csname newhelp\endcsname\athelp@
{Only certain combinations beginning with @ make sense to me.^^J%
I'll assume you wanted @@ for a printed @.}}{}
\@ifundefined{err@}{\def\err@{\@latexerr}}{}
\@ifundefined{default@tag}%
  {\def\default@tag{%
    \def\tag{\err@{\string\tag\space not allowed here}\@eha}}}
  {}
\@ifundefined{ex@}{\newdimen\ex@}{}
\@ifundefined{minaw@}{\newdimen\minaw@}{}
\@ifundefined{bigaw@}{\newdimen\bigaw@}{}
\newdimen\minCDarrowwidth
\newif\ifCD@
\let\ampersand@\relax
\def\CD{\catcode`\@\active
 \bgroup\relax\let\ampersand@&\iffalse}\fi
 \CD@true\vcenter\bgroup\Let@\restore@math@cr\default@tag
 \tabskip\z@skip\baselineskip20\ex@
 &\hfill$\m@th##$\hfill\crcr}
\def\endCD{\crcr\egroup\egroup\egroup}
\def\CD@check#1#2{\ifCD@\DN@{#2}\else
  \DN@{\err@{@\string#1 not
    allowed outside of the CD environment}\@eha}%
  \fi\next@}
\atdef@>#1>#2>{\ampersand@
  \ifCD@ \global\bigaw@\minCDarrowwidth \else \global\bigaw@\minaw@ \fi
  \setboxz@h{$\m@th\scriptstyle\;{#1}\;\;$}%
  \ifdim\wdz@>\bigaw@\global\bigaw@\wdz@\fi
  \@ifnotempty{#2}{\setbox@ne\hbox{$\m@th\scriptstyle\;{#2}\;\;$}%
    \ifdim\wd@ne>\bigaw@\global\bigaw@\wd@ne\fi}%
 \ifCD@\enskip\fi
   \mathrel{\mathop{\hbox to\bigaw@{\rightarrowfill@\displaystyle}}%
     \limits^{#1}\@ifnotempty{#2}{_{#2}}}%
 \ifCD@\enskip\fi \ampersand@}
\atdef@<#1<#2<{\ampersand@
  \ifCD@ \global\bigaw@\minCDarrowwidth \else \global\bigaw@\minaw@ \fi
  \setboxz@h{$\m@th\scriptstyle\;\;{#1}\;$}%
  \ifdim\wdz@>\bigaw@ \global\bigaw@\wdz@ \fi
  \@ifnotempty{#2}{\setbox@ne\hbox{$\m@th\scriptstyle\;\;{#2}\;$}%
    \ifdim\wd@ne>\bigaw@ \global\bigaw@\wd@ne \fi}%
  \ifCD@\enskip\fi
    \mathrel{\mathop{\hbox to\bigaw@{\leftarrowfill@\displaystyle}}%
      \limits^{#1}\@ifnotempty{#2}{_{#2}}}%
  \ifCD@\enskip\fi \ampersand@}
\begingroup \catcode`\~=\active \lccode`\~=`\@
\lowercase{%
  \global\atdef@)#1)#2){~>#1>#2>}
  \global\atdef@(#1(#2({~<#1<#2<}
}
\endgroup
\atdef@ A#1A#2A{\CD@check{A..A..A}{\llap{$\m@th\vcenter{\hbox
  {$\scriptstyle#1$}}$}\Big\uparrow
  \rlap{$\m@th\vcenter{\hbox{$\scriptstyle#2$}}$}&&}}
\atdef@ V#1V#2V{\CD@check{V..V..V}{\llap{$\m@th\vcenter{\hbox
  {$\scriptstyle#1$}}$}\Big\downarrow
  \rlap{$\m@th\vcenter{\hbox{$\scriptstyle#2$}}$}&&}}
\atdef@={\CD@check={&\enskip\mathrel
  {\vbox{\hrule\@width\minCDarrowwidth\vskip2\ex@\hrule\@width
  \minCDarrowwidth}}\enskip&}}
\atdef@|{\CD@check|{\Big\Vert&&}}
\atdef@\vert{\CD@check\vert{\Big\Vert&&}}
\atdef@.{\CD@check.{&&}}

\makeatletter

\newdimen\paperwidth
\newdimen\paperheight

\def\papersize#1#2{\let\p@persize\relax\paperwidth#1\paperheight#2}

\def\Afour{\papersize{210truemm}{297truemm}}

\let\p@persize\Afour

\let\onesidestyle\@twosidefalse
\let\twosidestyle\@twosidetrue

\def\margins{\@ifnextchar[{\@margins}{\@margins[\z@]}}

\def\@margins[#1]#2#3{
  \p@persize\dimen0 #3\dimen0 .5\dimen0\normalsize%
  \oddsidemargin-1truein\advance\oddsidemargin#2%
  \evensidemargin-1truein\advance\evensidemargin#2%
  \topmargin-1truein\advance\topmargin\dimen0\headsep\dimen0\footskip\dimen0%
  \textwidth\paperwidth\advance\textwidth-#2\advance\textwidth-#2%
  \textheight\paperheight\advance\textheight-#3\advance\textheight-#3%
  \headheight\baselineskip\advance\topmargin-.5\baselineskip%
  \advance\headsep-.5\baselineskip%
  \footheight\baselineskip
  \advance\textwidth-#1\advance\oddsidemargin#1
  \if@twoside\def\@themargin%
    {\ifodd\count\z@\oddsidemargin\else\evensidemargin\fi}\fi}

\def\headlinesep#1{\advance\topmargin\headsep\advance\topmargin -#1
  \advance\topmargin.5\baselineskip\headsep #1\advance\headsep-.5\baselineskip}

\def\footlinesep#1{\normalsize\footskip#1}

\def\headline{\if@twoside\let\n@xt\h@dlin@\else\let\n@xt\h@@dlin@\fi\n@xt}
  
\def\h@dlin@#1#2{%
  \def\@oddhead{%
    {{\leftskip\z@\rightskip\z@\noindent\normalsize#1}}}
  \def\@evenhead{%
    {{\leftskip\z@\rightskip\z@\noindent\normalsize#2}}}}

\def\h@@dlin@#1{%
  \def\@oddhead{{{\leftskip\z@\rightskip\z@\noindent\normalsize#1}}}}

\def\footline{\if@twoside\let\n@xt\f@tlin@\else\let\n@xt\f@@tlin@\fi\n@xt}

\def\f@tlin@#1#2{%
  \def\@oddfoot{%
    {{\leftskip\z@\rightskip\z@\noindent\normalsize#1}}}
  \def\@evenfoot{%
    {{\leftskip\z@\rightskip\z@\noindent\normalsize#2}}}}

\def\f@@tlin@#1{%
  \def\@oddfoot{{{\leftskip\z@\rightskip\z@\noindent\normalsize#1}}}}

\def\normalpage{\global\@specialpagefalse}

\makeatother

\makeatletter

\def\ft{\@ifnextchar[{\ft@s}{\ft@}}
\def\ft@{\ft@@@s[\f@size]}
\def\ft@s[{\@ifnextchar{a}{\ft@sz[}{\ft@@s[}}
\def\ft@@s[{\@ifnextchar{s}{\ft@sz[}{\ft@@@s[}}
\def\ft@@@s[#1]{\ft@sz[at #1pt]}
\def\ft@sz[#1]#2{\font\fonttemp=#2 #1\fonttemp\ignorespaces}

\makeatother

\makeatletter

\input{epsf.sty}

\def\showfig#1#2{\epsfbox{#2}}
\def\fig@#1#2{\leavevmode{\framebox{\figstyl@\strut{ #1 }}}}

\def\figstyle#1{\def\figstyl@{#1}}

\figstyle{\shape{n}\series{m}\selectfont\normalsize}

\def\showfigurestrue{\let\fig\showfig}
\def\showfiguresfalse{\let\fig\fig@}

\makeatother

\showfiguresfalse

\def\smallcircc{\mathop{\mkern3.5mu\hbox{\raise.58ex\hbox{\ft{lcircle10}a}}}}
\def\varemptyset{{\hbox{\raise.21ex\hbox{$\not$}}\mkern.15mu\mathrm{O}\mkern.15mu}}

\let\epsilon\varepsilon
\let\theta\vartheta
\let\phi\varphi

\let\emptyset\varemptyset

\documentstyle[12pt]{article}

\makeatletter

\let\Larg@\Large
\let\hug@\huge

\def\usepackage#1{\input{#1.sty}}

\input{geom.sty}

\let\Large\Larg@
\let\huge\hug@

\def\smallskip{\vskip\smallskipamount}
\def\medskip{\vskip\medskipamount}
\def\bigskip{\vskip\bigskipamount}

\def\mytrivlist{\parsep\parskip\@nmbrlistfalse
  \my@trivlist \labelwidth\z@ \leftmargin\z@
  \itemindent\z@ \def\makelabel##1{##1}}

\def\my@trivlist{\global\@newlisttrue \@outerparskip\parskip}

\def\end#1{\csname end#1\endcsname\@checkend{#1}%
  \expandafter\endgroup\if@endpe\@doendpe\fi
  \if@ignore \global\@ignorefalse \ignorespaces\fi}

\def\put{\@ifnextchar[{\@put}{\@@rput[\z@,\z@][r]}}
\def\@put[#1]{\@ifnextchar[{\@@put[#1]}{\@@@@@put[#1]}}
\def\@@put[#1][{\@ifnextchar{l}{\@@lput[#1][}{\@@@put[#1][}}
\def\@@@put[#1][{\@ifnextchar{c}{\@@cput[#1][}{\@@@@put[#1][}}
\def\@@@@put[#1][{\@ifnextchar{r}{\@@rput[#1][}{\relax}}
\def\@@@@@put[{\@ifnextchar{l}{\@@lput[\z@,\z@][}{\@@@@@@put[}}
\def\@@@@@@put[{\@ifnextchar{c}{\@@cput[\z@,\z@][}{\@@@@@@@put[}}
\def\@@@@@@@put[{\@ifnextchar{r}{\@@rput[\z@,\z@][}{\@@@@@@@@put[}}
\def\@@@@@@@@put[#1]{\@@rput[#1][r]}

\let\hm@d@\leavevmode

\long\def\@@lput[#1,#2][l]#3{\setbox0\hbox{#3}\hm@d@\raise#2\hbox to\z@{\dimen0 #1%
  \advance\dimen0-\wd0\kern\dimen0\dp0\z@\ht0\z@\wd0\z@\box0\hss}\ignorespaces}
\long\def\@@cput[#1,#2][c]#3{\setbox0\hbox{#3}\hm@d@\raise#2\hbox to\z@{\dimen0 #1%
  \advance\dimen0-.5\wd0\kern\dimen0\dp0\z@\ht0\z@\wd0\z@\box0\hss}\ignorespaces}
\long\def\@@rput[#1,#2][r]#3{\setbox0\hbox{\kern#1\raise#2\hbox{#3}}%
  \dp0\z@\ht0\z@\wd0\z@\hm@d@\box0\ignorespaces}

\def\flbox{\@ifnextchar[{\@flbox}{\@@rflbox[\z@,\z@][r]}}
\def\@flbox[#1]{\@ifnextchar[{\@@flbox[#1]}{\@@@@@flbox[#1]}}
\def\@@flbox[#1][{\@ifnextchar{l}{\@@lflbox[#1][}{\@@@flbox[#1][}}
\def\@@@flbox[#1][{\@ifnextchar{c}{\@@cflbox[#1][}{\@@@@flbox[#1][}}
\def\@@@@flbox[#1][{\@ifnextchar{r}{\@@rflbox[#1][}{\relax}}
\def\@@@@@flbox[{\@ifnextchar{l}{\@@lflbox[\z@,\z@][}{\@@@@@@flbox[}}
\def\@@@@@@flbox[{\@ifnextchar{c}{\@@cflbox[\z@,\z@][}{\@@@@@@@flbox[}}
\def\@@@@@@@flbox[{\@ifnextchar{r}{\@@rflbox[\z@,\z@][}{\@@@@@@@@flbox[}}
\def\@@@@@@@@flbox[#1]{\@@rflbox[#1][r]}
\long\def\@@lflbox[#1,#2][l]#3{\@@lput[#1,#2][l]{%
  \vtop{\leftskip\z@\parindent\z@\raggedleft\hm@d@#3}}}
\long\def\@@cflbox[#1,#2][c]#3{\@@cput[#1,#2][c]{%
  \vtop{\leftskip\z@\parindent\z@\raggedcenter\hm@d@#3}}}
\long\def\@@rflbox[#1,#2][r]#3{\@@rput[#1,#2][r]{%
  \vtop{\leftskip\z@\parindent\z@\raggedright\hm@d@#3}}}

\def\maketitle{\par
 \begingroup
 \def\thefootnote{\fnsymbol{footnote}}
 \def\@makefnmark{\hbox 
 to 0pt{$^{\@thefnmark}$\hss}} 
 \if@twocolumn 
 \twocolumn[\@maketitle] 
 \else 
 \global\@topnum\z@ \@maketitle \fi\thispagestyle{plain}\@thanks
 \endgroup
 \setcounter{footnote}{0}
 \let\maketitle\relax
 \let\@maketitle\relax
 \gdef\@thanks{}\gdef\@author{}\gdef\@title{}\let\thanks\relax}

\def\@maketitle{ 
 \null
 \vskip 2em \begin{center}
 {\LARGE \@title \par} \vskip 1.5em {\large \lineskip .5em
\begin{tabular}[t]{c}\@author 
 \end{tabular}\par} 
 \vskip 1em {\large \@date} \end{center}
 \par
 \vskip 1.5em}
 
\def\partbeforeskip#1{\def\p@rtbeforeskip{#1}}
\def\partstyle#1{\def\p@rtstyl@{#1}}
\def\partdot#1{\def\partd@t{#1}}
\def\partafterskip#1{\def\p@rtafterskip{#1}}
\def\partintrostyle#1{\def\partintr@styl@{#1}}
\def\partintrodot#1{\def\partintr@dot{#1}}
\long\def\partintrosep#1{\long\def\partintr@sep{#1}}
\def\partnewpagetrue{\def\p@rtnewp@ge{\newpage}}
\def\partnewpagefalse{\long\def\p@rtnewp@ge{\par}}

\partbeforeskip{4ex}
\partstyle{\centering\Large\bf}
\partdot{}
\partafterskip{3ex}
\partintrostyle{\large}
\partintrodot{}
\partintrosep{\par}
\partnewpagefalse

\def\partname{Part}
\def\part{\p@rtnewp@ge\addvspace\p@rtbeforeskip\@afterindentfalse\secdef\@part\@spart}

\def\@part[#1]#2{\ifnum \c@secnumdepth >-1\relax  
        \refstepcounter{part}                     
        \def\@tempa{\addcontentsline{toc}{part}}  %
        \expandafter\@tempa\expandafter{\thepart  
          \hspace{1em}#1}\else                    
        \addcontentsline{toc}{part}{#1}\fi        
   {\p@rtstyl@                       
    \ifnum \c@secnumdepth >-1\relax        
      {\partintr@styl@\partname\ \thepart  
       \partintr@dot}\partintr@sep\nobreak 
    \fi                                    
    #2\partd@t\markboth{}{}\par}
    \nobreak                       
    \vskip\p@rtafterskip           
   \@afterheading                  
    }                              

\def\@spart#1{{\p@rtcentering\p@rtstyl@                      
    #1\partd@t\par}                 
    \nobreak                        
    \vskip\p@rtafterskip            
    \@afterheading                  
  }                                 

\def\sectionbeforeskip#1{\def\s@ctbeforeskip{#1}}
\def\sectionstyle#1{\def\s@ctstyl@{#1}}
\def\sectiondot#1{\def\sectiond@t{#1}}
\def\sectionafterskip#1{\def\s@ctafterskip{#1}}
\def\sectionintrostyle#1{\def\sectionintr@styl@{#1}}
\def\sectionintro#1{\def\sectionintr@{#1}}
\def\sectionintrodot#1{\def\sectionintr@dot{#1}}
\def\sectionindenttrue{\def\s@ctind{\parindent}}
\def\sectionindentfalse{\def\s@ctind{\z@}}
\def\sectionafterindenttrue{\def\s@ct@ftind{+}}
\def\sectionafterindentfalse{\def\s@ct@ftind{-}}
\def\sectionafternewlinetrue{\def\s@ct@ftpar{+}}
\def\sectionafternewlinefalse{\def\s@ct@ftpar{-}}

\def\subsectionbeforeskip#1{\def\ss@ctbeforeskip{#1}}
\def\subsectionstyle#1{\def\ss@ctstyl@{#1}}
\def\subsectiondot#1{\def\subsectiond@t{#1}}
\def\subsectionafterskip#1{\def\ss@ctafterskip{#1}}
\def\subsectionintrostyle#1{\def\subsectionintr@styl@{#1}}
\def\subsectionintro#1{\def\subsectionintr@{#1}}
\def\subsectionintrodot#1{\def\subsectionintr@dot{#1}}
\def\subsectionindenttrue{\def\ss@ctind{\parindent}}
\def\subsectionindentfalse{\def\ss@ctind{\z@}}
\def\subsectionafterindenttrue{\def\ss@ct@ftind{+}}
\def\subsectionafterindentfalse{\def\ss@ct@ftind{-}}
\def\subsectionafternewlinetrue{\def\ss@ct@ftpar{+}}
\def\subsectionafternewlinefalse{\def\ss@ct@ftpar{-}}

\def\subsubsectionbeforeskip#1{\def\sss@ctbeforeskip{#1}}
\def\subsubsectionstyle#1{\def\sss@ctstyl@{#1}}
\def\subsubsectiondot#1{\def\subsubsectiond@t{#1}}
\def\subsubsectionafterskip#1{\def\sss@ctafterskip{#1}}
\def\subsubsectionintrostyle#1{\def\subsubsectionintr@styl@{#1}}
\def\subsubsectionintro#1{\def\subsubsectionintr@{#1}}
\def\subsubsectionintrodot#1{\def\subsubsectionintr@dot{#1}}
\def\subsubsectionindenttrue{\def\sss@ctind{\parindent}}
\def\subsubsectionindentfalse{\def\sss@ctind{\z@}}
\def\subsubsectionafterindenttrue{\def\sss@ct@ftind{+}}
\def\subsubsectionafterindentfalse{\def\sss@ct@ftind{-}}
\def\subsubsectionafternewlinetrue{\def\sss@ct@ftpar{+}}
\def\subsubsectionafternewlinefalse{\def\sss@ct@ftpar{-}}

\def\paragraphbeforeskip#1{\def\p@rbeforeskip{#1}}
\def\paragraphstyle#1{\def\p@rstyl@{#1}}
\def\paragraphdot#1{\def\paragraphd@t{#1}}
\def\paragraphafterskip#1{\def\p@rafterskip{#1}}
\def\paragraphintrostyle#1{\def\paragraphintr@styl@{#1}}
\def\paragraphintro#1{\def\paragraphintr@{#1}}
\def\paragraphintrodot#1{\def\paragraphintr@dot{#1}}
\def\paragraphindenttrue{\def\p@rind{\parindent}}
\def\paragraphindentfalse{\def\p@rind{\z@}}
\def\paragraphafterindenttrue{\def\p@r@ftind{+}}
\def\paragraphafterindentfalse{\def\p@r@ftind{-}}
\def\paragraphafternewlinetrue{\def\p@r@ftpar{+}}
\def\paragraphafternewlinefalse{\def\p@r@ftpar{-}}

\def\subparagraphbeforeskip#1{\def\sp@rbeforeskip{#1}}
\def\subparagraphstyle#1{\def\sp@rstyl@{#1}}
\def\subparagraphdot#1{\def\subparagraphd@t{#1}}
\def\subparagraphafterskip#1{\def\sp@rafterskip{#1}}
\def\subparagraphintrostyle#1{\def\subparagraphintr@styl@{#1}}
\def\subparagraphintro#1{\def\subparagraphintr@{#1}}
\def\subparagraphintrodot#1{\def\subparagraphintr@dot{#1}}
\def\subparagraphindenttrue{\def\sp@rind{\parindent}}
\def\subparagraphindentfalse{\def\sp@rind{\z@}}
\def\subparagraphafterindenttrue{\def\sp@r@ftind{+}}
\def\subparagraphafterindentfalse{\def\sp@r@ftind{-}}
\def\subparagraphafternewlinetrue{\def\sp@r@ftpar{+}}
\def\subparagraphafternewlinefalse{\def\sp@r@ftpar{-}}

\sectionbeforeskip{\bigskipamount}
\sectionstyle{\large\bf}
\sectiondot{}
\sectionafterskip{.5\bigskipamount}
\sectionintrostyle{}
\sectionintro{}
\sectionintrodot{.}
\sectionindentfalse
\sectionafterindenttrue
\sectionafternewlinetrue

\subsectionbeforeskip{.8\bigskipamount}
\subsectionstyle{\normalsize\bf}
\subsectiondot{}
\subsectionafterskip{.4\bigskipamount}
\subsectionintrostyle{}
\subsectionintro{}
\subsectionintrodot{.}
\subsectionindentfalse
\subsectionafterindenttrue
\subsectionafternewlinetrue

\subsubsectionbeforeskip{.6\bigskipamount}
\subsubsectionstyle{\normalsize\bf}
\subsubsectiondot{}
\subsubsectionafterskip{.3\bigskipamount}
\subsubsectionintrostyle{}
\subsubsectionintro{}
\subsubsectionintrodot{.}
\subsubsectionindentfalse
\subsubsectionafterindenttrue
\subsubsectionafternewlinetrue

\paragraphbeforeskip{.5\bigskipamount}
\paragraphstyle{\normalsize\bf}
\paragraphdot{.}
\paragraphafterskip{1.25ex}
\paragraphintrostyle{}
\paragraphintro{}
\paragraphintrodot{.}
\paragraphindentfalse
\paragraphafterindenttrue
\paragraphafternewlinefalse

\subparagraphbeforeskip{.5\bigskipamount}
\subparagraphstyle{\normalsize\bf}
\subparagraphdot{.}
\subparagraphafterskip{1.25ex}
\subparagraphintrostyle{}
\subparagraphintro{}
\subparagraphintrodot{.}
\subparagraphindenttrue
\subparagraphafterindenttrue
\subparagraphafternewlinefalse

\def\@startsection#1#2#3#4#5#6{
   \vskip\z@\@tempskipa #4\relax\@afterindenttrue
   \ifdim \@tempskipa <\z@ \@tempskipa -\@tempskipa \@afterindentfalse\fi
   \advance\@tempskipa by\presection
   \if@nobreak \everypar{}\else
     \addpenalty{\@secpenalty}\addvspace{\@tempskipa}%
     \allowbreak\vskip -\presection \fi \@ifstar
     {\@ssect{#1}{#2}{#3}{#4}{#5}{#6}}{\@dblarg{\@sect{#1}{#2}{#3}{#4}{#5}{#6}}}}

\def\@sect#1#2#3#4#5#6[#7]#8{\def\object@type{#1}%
   \ifnum #2>\c@secnumdepth\def\@svsec{}\def\@tempb{}%
      \else\refstepcounter{#1}\def\@svsec{{\csname #1intr@styl@\endcsname%
        {\csname #1intr@\endcsname}}\csname the#1\endcsname%
        \csname #1intr@dot\endcsname\kern1.25ex}%
        \edef\@tempb{\noexpand\numberline{\csname the#1\endcsname}}\fi%
   \@tempskipa #5\relax\def\@tempa{\addcontentsline{toc}{#1}}%
   \ifdim \@tempskipa>\z@%
      \begingroup #6\relax%
        \@hangfrom{\hskip #3\relax\@svsec}{\interlinepenalty \@M #8%
        \csname #1d@t\endcsname\par}%
      \endgroup%
      \csname #1mark\endcsname{#7}%
      \expandafter\@tempa\expandafter{\@tempb #7}%
      \ifautolabel\label*{#8}\fi%
   \else%
      \def\@svsechd{#6\hskip #3\relax%
         \@svsec #8\csname #1mark\endcsname {#7}%
         \expandafter\@tempa\expandafter{\@tempb #7}%
         \ifautolabel\label*{#8}\fi}\fi%
   \@xsect{#5}}

\def\@ssect#1#2#3#4#5#6#7{%
   \ifnum #2>\c@secnumdepth\def\@tempb{}\else \def\@tempb{\numberline{}}\fi%
     \@tempskipa #5\relax\def\@tempa{\addcontentsline{toc}{s#1}}%
     \ifdim \@tempskipa>\z@
        \begingroup #6\relax
           \@hangfrom{\hskip #3}{\interlinepenalty \@M #7%
           \csname #1d@t\endcsname\par}%
        \endgroup
        \csname s#1mark\endcsname{#7}%
        \ifstarredcontents\expandafter\@tempa\expandafter{\@tempb #7}\fi%
        \ifautolabel\label*{#7}\fi%
     \else%
        \def\@svsechd{#6\hskip #3\relax #7\csname s#1mark\endcsname {#7}%
        \ifautolabel\label*{#7}\fi}\fi
   \@xsect{#5}}

\def\section{\@startsection{section}{1}{\s@ctind}
  {\s@ct@ftind\s@ctbeforeskip}{\s@ct@ftpar\s@ctafterskip}{\s@ctstyl@}}
\def\subsection{\@startsection{subsection}{2}{\ss@ctind}
  {\ss@ct@ftind\ss@ctbeforeskip}{\ss@ct@ftpar\ss@ctafterskip}{\ss@ctstyl@}}
\def\subsubsection{\@startsection{subsubsection}{3}{\sss@ctind}
  {\sss@ct@ftind\sss@ctbeforeskip}{\sss@ct@ftpar\sss@ctafterskip}{\sss@ctstyl@}}
\def\paragraph{\@startsection{paragraph}{4}{\p@rind}
  {\p@r@ftind\p@rbeforeskip}{\p@r@ftpar\p@rafterskip}{\p@rstyl@}}
\def\subparagraph{\@startsection{subparagraph}{4}{\sp@rind}
  {\sp@r@ftind\sp@rbeforeskip}{\sp@r@ftpar\sp@rafterskip}{\sp@rstyl@}}

\def\statementabove#1{\def\th@bove{#1}}
\def\statementstyle#1{\def\thstyl@{#1}}
\def\statementbelow#1{\def\thb@low{#1}}
\def\statementindentfalse{\let\thind@nt\relax}
\def\statementindenttrue{\let\thind@nt\indent}

\def\statementintrostyle#1{\def\thintr@style{#1}}
\def\statementintrodot#1{\def\thintr@dot{#1}}
\def\statementintrosep#1{\def\thintr@sep{#1}}
\def\statementintrobrackets#1#2{\def\thintr@left{#1}\def\thintr@right{#2}}

\statementabove{\medskip}
\statementstyle{\sl}
\statementbelow{\medskip}
\statementindenttrue

\statementintrostyle{\normalshape\bf}
\statementintrodot{.}
\statementintrosep{\kern1.25ex}
\statementintrobrackets{(}{)}

\def\@thskip{\dimen0\lastskip\vskip-\dimen0%
  \th@bove\dimen1\lastskip\vskip-\dimen1%
  \ifdim\dimen0>\dimen1\else\dimen0\dimen1\fi\vskip\dimen0}

\long\def\@@newtheorem#1#2#3{%
  \newenvironment{#3}%
    {\def\object@type{#3}\par\@thskip#1%
     \@ifnextchar[{\@enva{#3}{\thstyl@{#2}}}{\@envb{#3}{\thstyl@{#2}}}}%
    {\end{#3@}}%
  \@ifnextchar[{\@othm{#3@}}{\@nnthm{#3}}}

\def\renewtheorem{\@ifnextchar[{\@renewtheorem}{\@renewtheorem[{}{}]}}

\long\def\@renewtheorem[#1]{\@@renewtheorem#1}

\long\def\@@renewtheorem#1#2#3{%
  \expandafter\let\csname#3@\endcsname\undefined
  \renewenvironment{#3}%
    {\def\object@type{#3}\par\@thskip%
     \@ifnextchar[{\@enva{#3}{\thstyl@#1{#2}}}{\@envb{#3}{\thstyl@#1{#2}}}}%
    {\end{#3@}}%
  \@ifnextchar[{\@othm{#3@}}{\@nnthm{#3}}}

\def\@begintheorem#1#2{\@opargbegintheorem{#1}{#2}{}}

\def\@opargbegintheorem#1#2#3{%
        \def\@tempx{#1}%
        \expandafter\let\expandafter\@tempy#2
        \def\@tempz{#3}%
        \mytrivlist\item[\thind@nt\hskip\labelsep%
        {\thintr@style%
          #1\if\@tempx\@empty\else\if\@tempy\relax\else\kern1ex\fi\fi#2%
          \ifx\@tempz\@empty%
            \if\@tempx\@empty\if\@tempy\relax%
            \else\thintr@dot\thintr@sep\fi\else\thintr@dot\thintr@sep\fi%
            \else%
            \if\@tempx\@empty\if\@tempy\relax\else\kern1ex\fi\else\kern1ex\fi%
           \thintr@left{#3}\thintr@right\thintr@dot\thintr@sep\fi}%
            \hskip-\labelsep]%
        \ifautolabel\label*{#3}\fi}

\def\@endtheorem{\strut\endtrivlist\thb@low}

\def\proofabove#1{\def\pf@bove{#1}}
\def\proofstyle#1{\def\pfstyl@{#1}}
\def\proofbelow#1{\def\pfb@low{#1}}
\def\proofindentfalse{\let\pfind@nt\relax}
\def\proofindenttrue{\let\pfind@nt\indent}

\def\proofintrostyle#1{\def\pfintr@style{#1}}
\def\proofintrodot#1{\def\pfintr@dot{#1}}
\def\proofintrosep#1{\def\pfintr@sep{#1}}
\def\proofintrobrackets#1#2{\def\pfintr@left{#1}\def\pfintr@right{#2}}

\proofabove{\medskip}
\proofstyle{}
\proofbelow{\medskip}
\proofindenttrue

\proofintrostyle{\sl}
\proofintrodot{.}
\proofintrosep{\kern1.25ex}
\proofintrobrackets{of\kern1ex}{}

\def\@pfskip{\dimen0\lastskip\vskip-\dimen0%
  \pf@bove\dimen1\lastskip\vskip-\dimen1%
  \ifdim\dimen0>\dimen1\else\dimen0\dimen1\fi\vskip\dimen0}

  {\@pfskip\mytrivlist\item[\pfind@nt]\@ifnextchar[{\pro@f}{\pro@f[\prooftag]}}
  {\ifvoid\provedbox\else\hproved\fi\endtrivlist\pfb@low}

\def\pro@f[#1]{\setbox\provedbox\hbox{\provedboxcontents{#1}}\proofintro{#1}}

\def\proofintro#1{\expandafter\def\expandafter\@tempa\expandafter{#1}%
  {\pfintr@style{Proof\ifx\@tempa\empty\else\kern1ex\pfintr@left{#1}%
  \pfintr@right\fi}\pfintr@dot\pfintr@sep}\pfstyl@\ignorespaces}

\def\provedmark#1{\def\prm@rk{#1}}
\def\provedsep#1{\def\prs@p{#1}}

\provedmark{$\square$}
\provedsep{\kern1.25ex}

\def\provedtexttrue{\def\prb@x##1{\fbox{\small##1}}}
\def\provedtextfalse{\def\prb@x##1{\prm@rk}}
\def\provedmarkrighttrue{\let\prhf@l\hfill}
\def\provedmarkrightfalse{\let\prhf@l\relax}

\provedtextfalse
\provedmarkrightfalse

\def\provedboxcontents#1{\expandafter\def\expandafter\@tempa\expandafter{#1}%
  \ifx\@tempa\empty\prm@rk\else\prb@x{#1}\fi}

\def\proved{\ifmmode\eqno{\box\provedbox}\else\hproved\fi}

\def\hproved{\unskip\nobreak\prhf@l\penalty50\prs@p\hbox{}\nobreak\prhf@l
  \box\provedbox{\finalhyphendemerits=0\par}}

\def\captionstyle#1{\def\c@ptstyl@{#1}}
\def\captionintrostyle#1{\def\c@pintr@style{#1}}
\def\captionintrodot#1{\def\c@pintr@dot{#1}}
\def\captionintrosep#1{\def\c@pintr@sep{#1}}

\captionstyle{\small\sf}
\captionintrostyle{\bf}
\captionintrodot{.}
\captionintrosep{\hskip1.25ex}

\long\def\@makecaption#1#2{%
    \vskip\captionskip
    \setbox\@tempboxa\hbox{%
      \ifproofing\@ifundefined{the@label}{}
        {\hbox to 0pt{\vbox to 0pt{\vss\hbox{\tiny\the@label}\bigskip}\hss}}\fi
      \c@ptstyl@{\c@pintr@style #1\c@pintr@dot}\ignorespaces #2}%
    \@captionwidth=\hsize \advance\@captionwidth-2\@captionmargin
    \ifdim \wd\@tempboxa >\@captionwidth {%
        \rightskip=\@captionmargin\leftskip=\@captionmargin
        \unhbox\@tempboxa\par}%
      \else
        \hbox to\hsize{\hfil\box\@tempboxa\hfil}%
    \fi}

\def\end@Float#1{%
  \expandafter\caption\expandafter[\the@title]{%
   {\c@pintr@style%
   \ifx\the@caption\empty\ifx\the@title\empty\else\c@pintr@sep\fi\else\c@pintr@sep\fi
    \the@title\ifx\the@caption\empty\else\ifx\the@title\empty\else
    \c@pintr@dot\c@pintr@sep\fi\fi}%
   \ignorespaces\the@caption%
  \expandafter\label\expandafter*\expandafter{\the@label}}%
  \end{#1}}

\@definecounter{bibenumi}

\def\thebibliography#1{\section*{\refname}%
 \list{[\arabic{bibenumi}]}{\settowidth\labelwidth{[#1]}%
 \leftmargin\labelwidth\advance\leftmargin\labelsep\usecounter{bibenumi}}%
 \def\newblock{\hskip .11em plus .33em minus .07em}%
 \sloppy\clubpenalty4000\widowpenalty4000\sfcode`\.=1000\relax}

\makeatother

\proofingfalse
\autolabelfalse

\showfiguresfalse

\newtheorem{stat}{\statname}  \unnumbered{stat}
\newenvironment{statement}[1]{\def\statname{#1}\begin{stat}}{\end{stat}}

\newtheorem{nstat}{\nstatname}[section]

\newtheorem{question}[nstat]{Question}

\newtheorem[{\ns}{}]{remark}[nstat]{Remark}

\papersize{215truemm}{275truemm}
\margins{3.25cm}{1.75cm}
\footlinesep{1cm}
\headline{\hfill}
\footline{\small\hfill--\kern1ex\thepage\kern1ex--\hfill}
\showfigurestrue

\sectionbeforeskip{1.5\bigskipamount}
\sectionstyle{\centering\normalsize\bf}
\sectionafterskip{\bigskipamount}
\sectionafterindenttrue

\paragraphbeforeskip{\bigskipamount}
\paragraphstyle{\centering\normalsize\sl}
\paragraphafterskip{.5\bigskipamount}
\paragraphafternewlinetrue
\paragraphafterindenttrue
\paragraphdot{}

\statementintrostyle{\sc}

\renewtheorem{question}{Question}

\captionstyle{\small}
\captionintrostyle{\sc}

\newcommand{\B}{B} 
\renewcommand{\S}{S}

\newcommand{\Cl}{\mathop{\mathrm{Cl}}\nolimits} 
\newcommand{\Int}{\mathop{\mathrm{Int}}\nolimits} 
\newcommand{\Bd}{\mathop{\mathrm{Bd}}\nolimits}

\renewcommand{\Im}{\mathop{\mathrm{Im}}\nolimits} 

\def\(#1){({\sl #1\/})}

\let\emptyset\varemptyset

\font\ftt cmtt10 at 11pt
\font\bfs cmbxsl10 at 12pt

\begin{document}


\title{\large\bf A UNIVERSAL RIBBON SURFACE IN $\B^4$%
       \label{Version 1.0 / \today}}
\author{\sc\normalsize R. Piergallini\\
\sl\normalsize Dipartimento di Matematica e Informatica\\[-3pt]
\sl\normalsize Universit\`a di Camerino -- Italia\\
\ftt riccardo.piergallini@unicam.it
\and 
\sc\normalsize D. Zuddas\\
\sl\normalsize Scuola Normale Superiore di Pisa -- Italia\\
\ftt d.zuddas@sns.it}
\date{}

\maketitle

\begin{abstract}
\baselineskip13.5pt
\smallskip
\noindent
We construct an orientable ribbon surface $F \subset B^4$, which
is universal in the following sense: any orientable 4-manifold
$M \cong B^4 \cup \text{1-handles} \cup \text{2-handles}$
can\break be represented as a cover of $B^4$ branched over $F$.



\end{abstract}


\section*{Introduction}

In the early seventies H.M. Hilden, U. Hirsch and J.M. Montesinos independently
proved that any closed orientable 3-manifold can be represented as a 3-fold simple
covering of $S^3$ branched over a knot (cf. \cite{Hi74}, \cite{H74} and
\cite{M74}). 

Ten years later, W. Thurston constructed the first universal link. He called a link
$L \subset S^3$ universal iff for any closed orientable 3-manifold there exists an
$n$-fold (in general non-simple) covering $M \to S^3$ branched over $L$.
Subsequently, other universal links and knots were constructed by H.M. Hilden, M.T.
Lozano, J.M. Montesinos and W.C. Whitten. The basic idea of these constructions is
the following: symmetrize the branching links given by the Hilden-Hirsch-Montesinos
representation theorem, making them sublinks of the preimage of a fixed link with
respect to a fixed branched covering $S^3 \to S^3$. (cf. \cite{T82}, \cite{HLM83},
\cite{HLM85}, \cite{HLM85'} and \cite{HLMW87}).

More recently, M. Iori and R. Piergallini obtained a representation theorem of
closed orientable smooth 4-manifolds as 5-fold simple covering of $S^4$ branched
over a smooth surface (cf. \cite{P95} and \cite{IP02}). Then, it makes sense to look
for a universal surface in $S^4$, satisfying a universal property analogous to that
one of a universal link in the 3-dimensional case. But unfortunately, the
symmetrization technique used for branching links in $S^3$ seems hardly to be
directly adaptable to branching surfaces in $S^4$.

\break

In this paper, we show how certain ribbon branching surfaces in $B^4$ can be
symmetrized, in order to get a universal orientable ribbon surface, for representing
any compact bounded orientable 4-manifold $M \cong B^4 \cup \text{1-handles} \cup
\text{2-handles}$ as a branched cover of $B^4$. Such 4-manifolds turn out to be
relevant for the presentation of all the closed orientable smooth 4-manifolds,
making no difference how 3- and 4-handles are attached to them (cf. \cite{K89}).
Hence, our result could be also useful in constructing a universal surface in $S^4$.
Namely, we prove the following theorem.

\begin{statement}{Theorem}
There exists an orientable ribbon surface $F \subset B^4$, such that any compact
orientable 4-manifold $M \cong B^4 \cup \text{1-handles} \cup \text{2-handles}$ is a
cover of $B^4$ branched over $F$. 
\end{statement}

We recall that a smooth compact surface $F \subset B^4$ with $\Bd F \subset S^3$ is
called a {\sl ribbon surface} if the Euclidean norm restricts to a Morse function on
$F$ with no local maxima in $\Int F$. Assuming $F \subset R^4_+ \subset R^4_+ \cup
\{\infty\} \cong B^4$, this property is topologically equivalent to the fact that
the fourth Cartesian coordinate restricts to a Morse function on $F$ with no local
minima in $\Int F$. Such a surface $F \subset R^4_+$ can be isotoped to make its
orthogonal projection $F_0 \subset R^3$ a self-transversal immersed surface, whose
double points form disjoint arcs as in Figure \ref{ribbon/fig}. Here, as well as in
the following figures, we shade the surface roughly according to the fourth
coordinate.

\begin{Figure}[htb]{ribbon/fig}{}{}\vskip1mm
\centerline{\fig{Ribbon self-intersection}{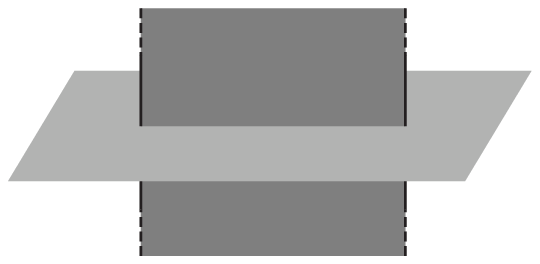}}
\end{Figure}

We will refer to $F_0$ as a 3-dimensional {\sl diagram} of $F$. It is worth
observing that any immersed compact surface $F_0 \subset R^3$ with no closed
components, all self-intersections of which are as above, is the diagram of a ribbon
surface $F$ uniquely determined up to isotopy. In fact, $F$ can be obtained by
pushing $\Int F_0$ inside $\Int B^4$, in such a way that all the self-intersections
disappear.

We also recall that a smooth map $p:M \to B^4$ is called a {\sl branched covering\/}
with {\sl branching surface} $F \subset B^4$, if the restriction $p_|: M - p^{-1}(F)
\to B^4 - F$ is a regular (as a smooth map) ordinary covering of finite degree $d$.
Assuming $F$ minimal with respect to this property, we have $F = p(S)$ where $S
\subset M$ is the {\sl singular surface} of $p$. For each $x \in S$, there exists a
neighborhood $U$ of $x$ in $M$, such that the restriction $p_|:U \to p(U)$ is
topological equivalent to the complex map defined by $(z_1,z_2) \mapsto
(z_1^{d_x},z_2)$ (with $\Im z_2 \geq 0$ if $x \in \Bd M$), where $d_x$ is a positive
integer called {\sl local degree} or {\sl branching index} of $p$ at $x$. We say
that $p$ is {\sl simple} if the restriction $p_|:S \to F$ is injective and $d_x = 2$
for any $x \in S$.

Since $p$ is uniquely determined (up to diffeomorphism) by its restriction over $B^4
- F$, we can describe it in terms of its {\sl monodromy} $\pi_1(B^4-F) \to
\Sigma_{d}$ (defined up to conjugation in $\Sigma_{d}$, depending on the numbering
of the sheets), that is by giving the monodromies of any set of meridian loops around
$S$ generating $\pi_1(B^4 - F)$. Usually, this is done by labelling a 3-dimensional
diagram of $F$ with the permutations corresponding to the standard Wirtinger's
generators of $\pi_1(B^4 - F)$, which are assumed to cross the diagram from back to
front. From this perspective, $p$ is simple if and only if all the labels are
transpositions.

\medskip

The paper is entirely devoted to prove the theorem above. In particular, the
symmetrization procedure is described in Section \ref{surface/sec} and the universal
surface $F$ is depicted in Figure \ref{univ/fig}. Sections \ref{covering/sec} and
\ref{moves/sec} are respectively aimed to show that any 4-manifold $M$ as
in the statement is a simple covering of $B^4$ branched over a suitable ribbon
surface and to introduce the covering moves needed for symmetrizing such a ribbon
branching surface.


\section{Some covering moves\label{moves/sec}}

By a {\sl covering move}, we mean any modification on a labelled surface determining
a branched covering $p:M \to B^4$, that preserves the covering manifold $M$ up to
diffeomorphism. All the covering moves considered in this paper are {\sl local},
that is the modification takes place inside a cell and can be performed whatever is
the rest of labelled branching surface outside. In the figures describing these
moves, we will draw only the part of the labelled branching surface inside the
relevant cell, assuming everything else to be fixed.

Of course, the notion of covering move makes sense for coverings between PL
manifolds of any dimension $m$ branched over arbitrary $(m-2)$-dimensional
subcomplexes of the range. Before of defining our moves, we roughly state two very
general equivalence principles in this broader context and discuss some applications
to our specific situation. Several special cases of these principles already appeared
in the literature and we can think of them as belonging to the ``folklore'' of
branched coverings.

\begin{statement}{Disjoint monodromies crossing}
Subcomplexes of the branching set of a covering that are labelled with disjoint
permutations can be isotoped independently from each other without changing the
covering manifold.
\end{statement}

The reason why this principle holds is quite simple. Namely, being the labelling of
the considered subcomplexes disjoint, the sheets non-trivially involved by them do
not interact, at least locally over the region where the isotopy take place. Hence,
relative position of such subcomplexes is not relevant in determining the covering
manifold.

\begin{Figure}[htb]{crossing/fig}{}{}
\centerline{\fig{Crossing change}{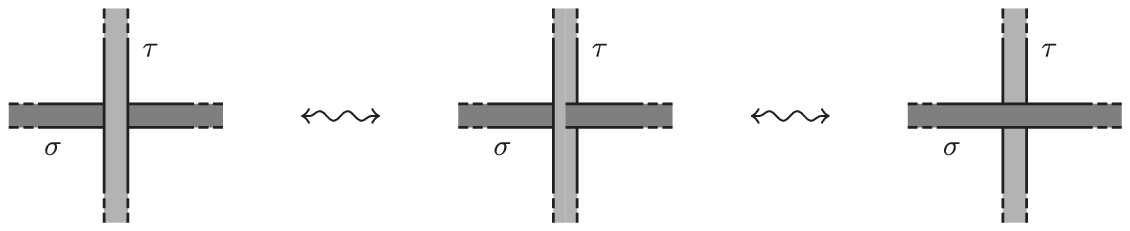}}
\end{Figure}

In particular, this principle allows crossing changes in diagrams when the involved
monodromies are disjoint. For example, this is the case of one of the well-known
Montesinos moves (cf. \cite{M85}, \cite{P95}, \cite{Ap03} or \cite{BP03}) for simple
coverings of $S^3$ branched over a link. Such crossing change has already been used
in the construction of uni\-versal links (cf. \cite{HLM85} and \cite{HLMW87}).
In the same spirit, we specialize the above principle in our 4-dimensional context,
by considering the crossing change move described in Figure \ref{crossing/fig},
where $\sigma, \tau \in \Sigma_d$ are arbitrary disjoint permutations. 

It is worth observing that, abandoning transversality, the disjoint monodromies
crossing principle also gives the special case of the next one when the $\sigma_i$'s
are disjoint and $L$ is empty.

\begin{statement}{Coherent monodromies merging}
Let $p:M \to N$ be any branched covering with branching set $B_p$ and let $\pi: E \to
K$ be a connected disk bundle imbedded in $N$, in such a way that: 1)~there exists a
(possibly empty) subcomplex $L \subset K$ for which $B_p \cap \pi^{-1}(L) = L$ and
the restriction of $\pi$ to $B_p \cap \pi^{-1}(K-L)$ is an unbranched covering of
$K-L$; 2)~the monodromies $\sigma_1, \dots, \sigma_n$ relative to a fundamental
system $\omega_1, \dots, \omega_n$ for the restriction of $p$ over a given disk
$D = \pi^{-1}(x)$, with $x \in K - L$, are coherent in the sense that $p^{-1}(D)$ is
a disjoint union of disks. Then, by contracting the bundle $E$ fiberwise to $K$, we
get a new branched covering $p':M \to N$, whose branching set $B_{p'}$ is equivalent
to $B_p$, except for the replacement of $B_p \cap \pi^{-1}(K-L)$ by $K-L$, with the
labelling uniquely defined by letting the monodromy of the meridian $\omega =
\omega_1 \dots \omega_n$ be $\sigma = \sigma_1 \dots \sigma_n$.
\end{statement}

We remark that, by connection and property 1), the coherence condition required in
2) actually holds for any $x \in K$. Then, we can prove that $p$ and $p'$ have the
same covering manifold, by a straightforward fiberwise application of the Alexander's
trick to the components of the bundle $\pi \circ p: p^{-1}(E) \to K$. A coherence
criterion can be immediately derived from Section 1 of \cite{MP01}.

We will mainly apply the merging principle to ``parallel'' components of the
branching surface with coherent monodromies, in order to control the number of such
components (cf. the below discussion of stabilization and Figures \ref{move2a/fig},
\ref{move3a/fig}, \ref{move3b/fig}, \ref{move4a/fig}).

However, this principle originated from a classical perturbation argument in
algebraic geometry and appeared in the literature as a way to deform non-simple
coverings between surfaces into simple ones, by going in the opposite direction from
$p'$ to $p$ (cf. \cite{BE79}). In dimension 3, it can be used in this direction, not
only for achieving simplicity (cf. \cite{BP03} or \cite{Hr03}), but also for removing
singularities from the branching set (cf. \cite{BP03}). Moreover, it has been used
in the construction of universal links, for controlling the branching indices (cf.
\cite{HLMW87}).

\begin{Figure}[htb]{merge/fig}{}{}
\centerline{\fig{Merging}{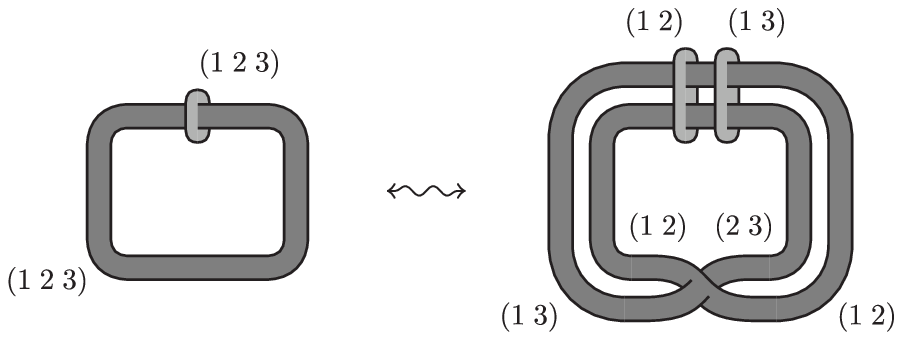}}
\end{Figure}

Figure \ref{merge/fig} shows an example of application of the merging principle to
coverings of $B^4$ branched over ribbon surfaces. Here, the absolute version (with
$L = \emptyset$) of the principle is applied in turn to both the components of the
branching surface on the left side (letting $K$ be a component and $\pi:E \to K$ be
its normal bundle). There is no obstruction to generalize this example, to show that
any covering of $B^4$ branched over a ribbon surface can be deformed into a simple
one. For applications of the relative version of the principle (in both
directions) see Figures \ref{move2a/fig} and \ref{move3b/fig}.

\medskip

Now, we pass to define our moves on labelled ribbon surfaces representing branched
coverings of $B^4$. Concerning the assumptions on the monodromies, the definitions
are given on a level of generality which is not the highest possible, but is still
higher than needed for our present purposes. We made this choice because such moves
are interesting in their own right. In the next section we will use only
stabilization and Moves 3 and 4. Moves 1 and 2 are used here to get the other ones.
Let us start with some considerations about the well-known notion of stabilization. 

\medskip

{\bfs Stabilization.} The basic version consists in the addition of an extra
trivial sheet, the $(d+1)$-th one, to a given $d$-fold branched covering. In terms of
branching surface, this means to add a separate trivial disk with label
$(i\ d{+}1)$, where $1 \leq i \leq d$. Now, we can iterate this process $l$ times,
by adding $l$ trivial disks with labels $(i_1\;d{+}1), \dots, (i_l\;d{+}l)$, where $1
\leq i_j \leq d + j - 1$. Of course, we can assume the disks to be parallel and
it is easy to realize that their monodromies are coherent, whatever $i_j$'s we
choose. Hence, we can merge all the disks into one. In particular, if all the $i_j$'s
are distinct, the label of this disk is given by the product of $l$ disjoint
transposition $(i_1\;d{+}1) \dots (i_l\;d{+}l)$. We will refer to the addition of
such a labelled disk as the {\sl multi-stabilization} involving the sheets $i_1,
\dots, i_l$.
\smallskip

\medskip

{\bfs Move 1.}
This move is described in Figure \ref{move1/fig}, where $j_1, \dots, j_l$ and 
$k_1, \dots, k_l$ are assumed to be all distinct (cf. \cite{HLMW87} and%
\begin{Figure}[htb]{move1/fig}{}{}\vskip1mm
\centerline{\fig{Move 1}{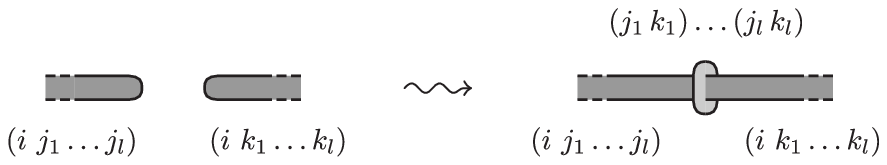}}
\end{Figure}
\cite{IP02} for the case of $l = 1$). It can be obtained by a straightforward
application of the main technique of \cite{IP02}, that is by extending the covering 
in the left side of the figure to certain cancelling 1- and 2- handles added to $M$
and $B^4$, in such a way that the branching surface becomes as in the right side.

Namely, we add to $B^4$ a 1-handle $H^1$, connecting two small 3-balls around the
tips of the tongues in the left side of the figure, and then a 2-handle $H^2$
complementary to $H^1$, whose attaching loop $\lambda$ meets $B^4$ along a horizontal
line avoiding the tongues. The covering instructions can be extended to these
handles, by assigning to $\lambda$ the monodromy $(j_1\;k_1) \dots (j_l\;k_l)$ and
by completing the branching surface with the cocore disk of $H^2$ labelled with the
same monodromy of $\lambda$. After cancelling $H^1$ and $H^2$, the new branching
surface and monodromy look like in the right side of Figure \ref{move1/fig}. We
leave to the reader to check that, in the new covering manifold, there are $d - l$
1-handles over $H^1$ and the same number of 2-handle over $H^2$ and that they cancel
(non-trivially) to give back $M$ again.

We remark that Move 1 could also be derived from the special case when $l = 1$, with
an inductive argument analogous to the one used below for Move 3.

\medskip

{\bfs Move 2.}
Our second move is given by Figure \ref{move2/fig}. Here, the $\sigma$ in the left
side is any permutation in $\Sigma_d$, while the $\sigma$ in the right side is the%
\begin{Figure}[htb]{move2/fig}{}{}\vskip1.5mm
\centerline{\fig{Move 2}{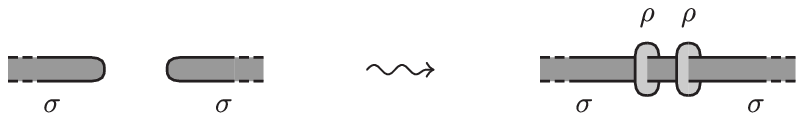}}
\end{Figure}
same permutation thought in $\Sigma_{d'}$, for a certain $d' > d$, and $\rho \in
\Sigma_{d'}$ is a product of disjoint transpositions which depends on $\sigma$.
Differently from the previous one, this move changes the degree of the covering. In
fact, we can transform the left side of Figure \ref{move2/fig} into the right one,
by performing a suitable multi-stabilization followed by a generalized version of
our first move. Let $\sigma = \gamma_1 \dots \gamma_h$ a decomposition of the
given permutation $\sigma$ into disjoint cycles. For the sake of exposition, we
proceed by induction on $h$.

If $h = 1$ we can write $\sigma = (i\ j_1\,\dots\,j_l)$. In this case, we
perform on the covering represented by the diagram on the left side of Figure
\ref{move2/fig} a multi-stabilization involving the sheets $j_1, \dots, j_l$. As a
result, one trivial disk with label $\rho = (j_1\;d{+}1) \dots (j_l\;d{+}l)$
appears in the diagram. We stretch one of the two tongues to pass through such
disk, so that its monodromy beyond it becomes $\sigma^\rho = \rho^{-1}\sigma\rho =
(i\ d{+}1\,\dots\,d{+}l)$. At this point, Move 1 immediately gives the diagram on the
right side of the figure.

The case of $h > 1$ can be reduced to the inductive hypothesis by means of crossing
changes and merging principle, as shown in Figure \ref{move2a/fig}. Here, $\sigma' =
\gamma_1 \dots \gamma_{h-1}$, $\sigma'' = \gamma_h$, $\rho'$ (resp. $\rho''$) is
the product of disjoint transpositions resulting from applications of Move 1 to the
tongues labelled with $\sigma'$ (resp. $\sigma''$), and $\rho = \rho' \rho''$.%
\begin{Figure}[htb]{move2a/fig}{}{}\vskip1.5mm
\centerline{\fig{Move 2 (proof)}{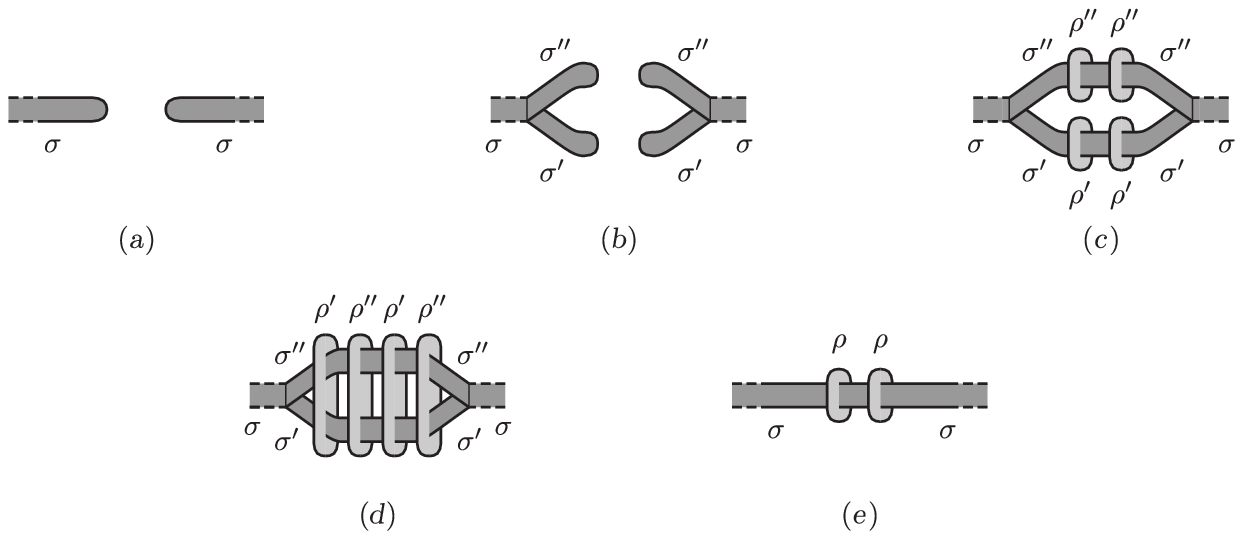}}
\end{Figure}
Starting from \(a), we apply in sequence: merging principle to get \(b); inductive
hypothesis to get \(c); crossing changes to get \(d); merging principle again to get
\(e).

\medskip

{\bfs Move 3.} Our third move is the one of Figure \ref{move3/fig}, where the
permutations $\sigma$ and $\rho$, as well as covering degrees, are the same of
Figure \ref{move2/fig}.%
\begin{Figure}[b]{move3/fig}{}{}
\centerline{\fig{Move 3}{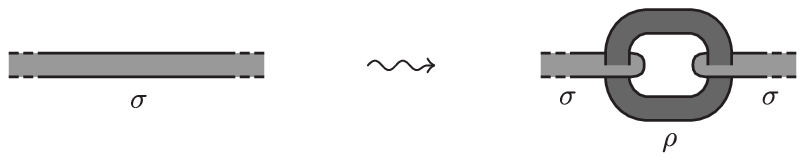}}
\end{Figure}
We can limit ourselves to consider the case when $\sigma$ is a cycle, since the
general case can be derived by induction on the length of a cyclic decomposition
of $\sigma$, with the same argument used for Move 2 (think of Figure
\ref{move3b/fig} below, as if it were labelled analogously to Figure
\ref{move2a/fig}). 

So, we assume $\sigma = (i\ j_1\,\dots\,j_l)$ and proceed by induction on $l$.
Figure \ref{move3a/fig} shows how to deal with the case of $l = 1$, when
$\sigma = (i\ j_1)$ and $\rho = (j_1\;d{+}1)$.\break We observe that diagrams \(c)
and \(d) represent isotopic surfaces, and same holds for diagrams \(e), \(f) and
\(g). Moreover: \(b) is a stabilization of \(a); \(c) and \(i) are obtained from
the previous diagrams by Move 2 and its inverse; \(e) and \(h) by crossing
changes. The inductive step is described in Figure \ref{move3b/fig}. Here, the
sequence of operations needed to get the various diagrams is the same of Figure
\ref{move2a/fig}.%
\begin{Figure}[htb]{move3a/fig}{}{}
\centerline{\fig{Move 3 (proof 1)}{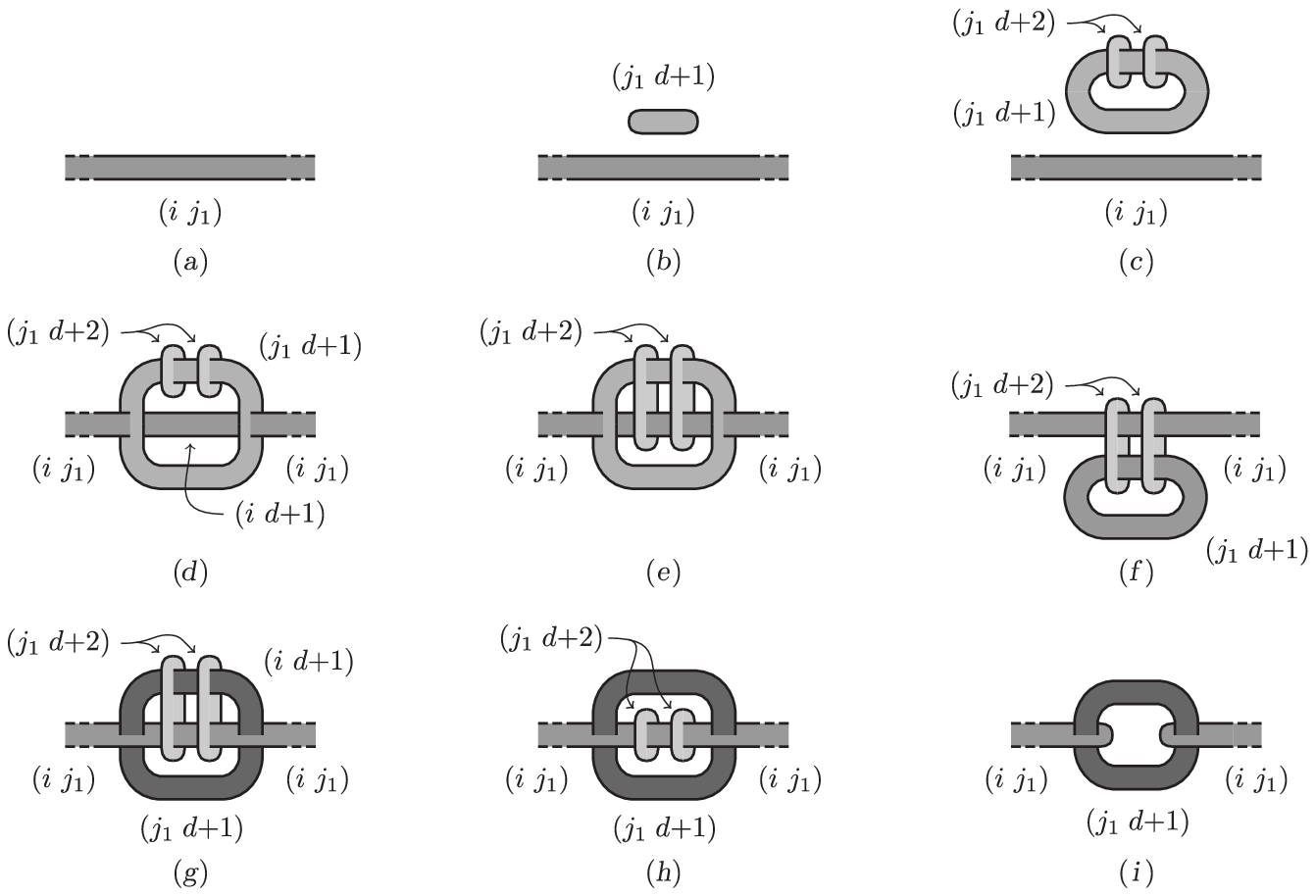}}\vskip-1mm
\end{Figure}
\begin{Figure}[htb]{move3b/fig}{}{}\vskip-4mm
\centerline{\fig{Move 3 (proof 2)}{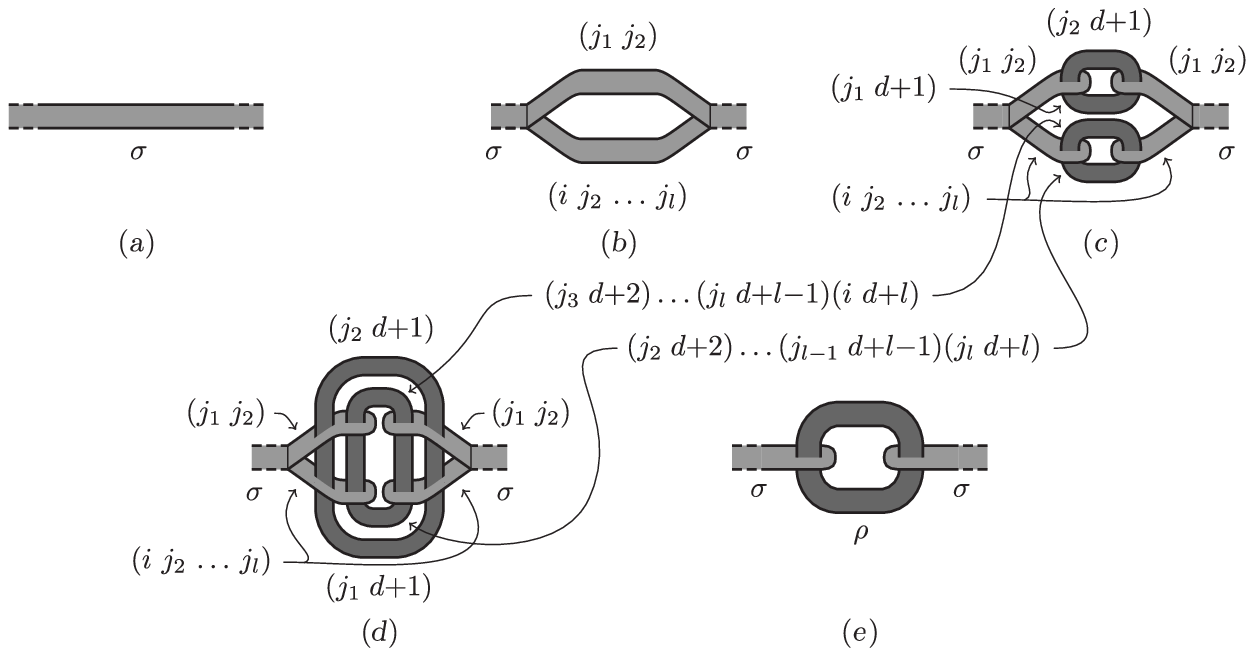}}\vskip-1mm
\end{Figure}


\break

{\bfs Move 4.} Differently from the previous ones, this move is defined only for
simple monodromies, but it does not preserve simplicity. It is depicted in Figure
\ref{move4/fig}, where $\tau_1, \tau_2 \in \Sigma_d$ are arbitrary distinct
transpositions and $\tau_3 = \tau_1^{\tau_2} = \tau_2^{-1} \tau_1 \tau_2$, while each
$\rho_j$ is a product of two disjoint transpositions which depend on the $\tau_i$'s.

\begin{Figure}[htb]{move4/fig}{}{}\vskip1mm
\centerline{\fig{Move 4}{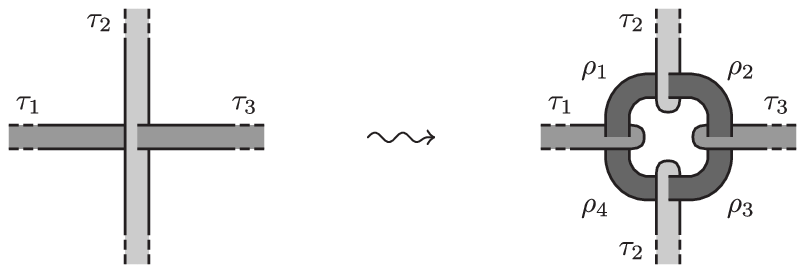}}
\end{Figure}

The following Figure \ref{move4a/fig} tells us why this is a true covering move
if $\tau_1$ and $\tau_2$ are not disjoint. Here, \(b) and \(c) are obtained by Move
3 (followed by isotopy in the former step), \(d) by crossing change, and \(e)
by merging principle. We leave to the reader to adapt the monodromies of Figure
\ref{move4a/fig} to the easier case of $\tau_1$ and
$\tau_2$ disjoint.

\begin{Figure}[htb]{move4a/fig}{}{}
\centerline{\fig{Move 4 (proof)}{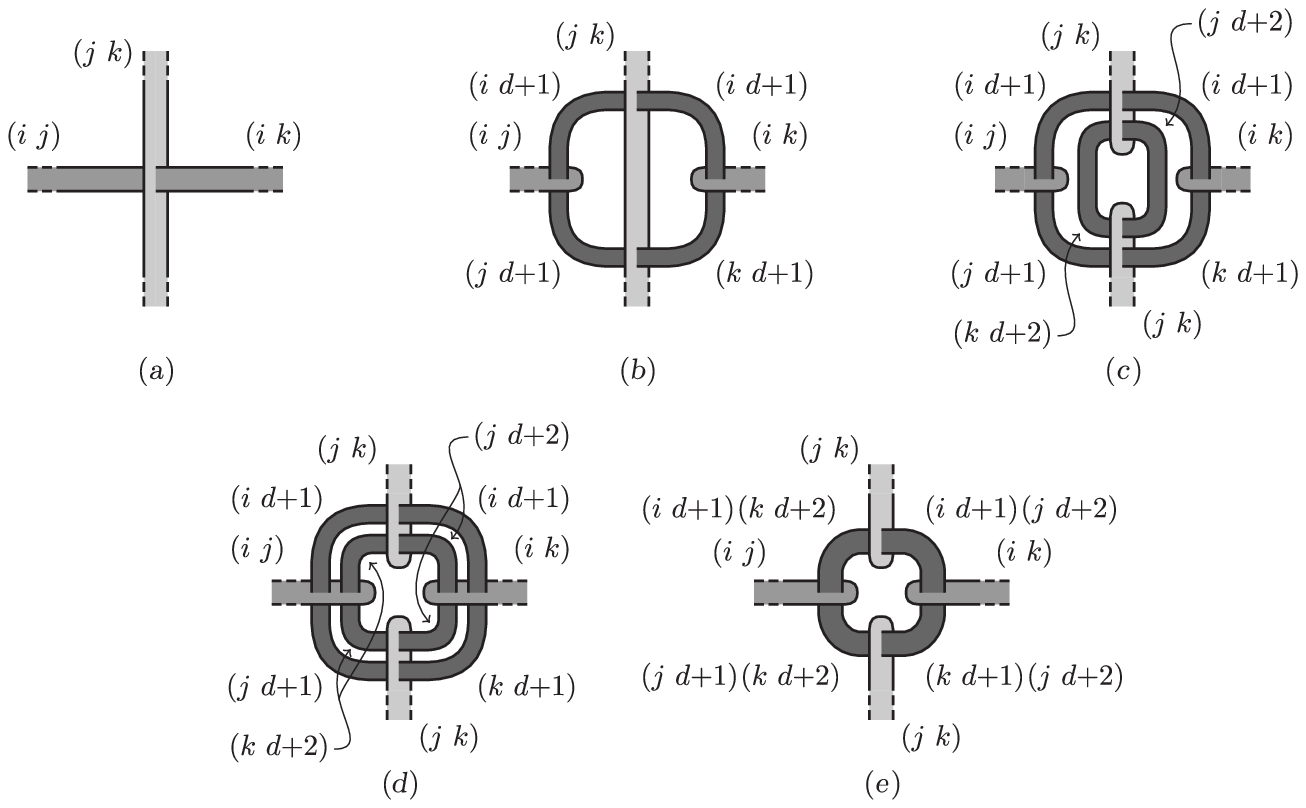}}
\end{Figure}


\section{Special covering presentations\label{covering/sec}}

Given a compact bounded orientable 4-manifold $M$ as in the statement of our
theorem, that is $M \cong B^4 \cup \text{1-handles} \cup \text{2-handles}$, we want
to present it as a simple covering of $B^4$ branched over a suitable orientable
ribbon surface.

By Montesinos \cite{M78}, we know that $M$ is a 3-fold simple cover of $B^4$
branched over a possibly non-orientable ribbon surface $F \subset B^4$. A variation 
of the Montesinos's argument actually shows that $F$ can be chosen orientable.
Alternatively, we can think of $M$ as a topological Lefschetz fibration over $B^2$
and represent it as 3-fold simple cover of $B^4$ branched over a braided surface
(cf. Remark 3 of \cite{LP01}).

However, we will construct a special covering presentation of $M$ by a technique
similar to one used in \cite{LP01}. This choice, renouncing to control the degree
of the covering, which is not relevant in this context, will eventually allows us to
get a simpler universal surface. Nevertheless, it is worth observing that the
symmetrization process described in the next section could be arranged to work
starting from a generic ribbon branched surface.

\medskip

Let the generic Kirby diagram of Figure \ref{diag1/fig} represent $M$. Here,
as well as in Figure \ref{diag2/fig}, the framings are assumed to coincide with the
blackboard ones outside the box.

\begin{Figure}[htb]{diag1/fig}{}{}\vskip-1mm
\centerline{\fig{Generic Kirby diagram}{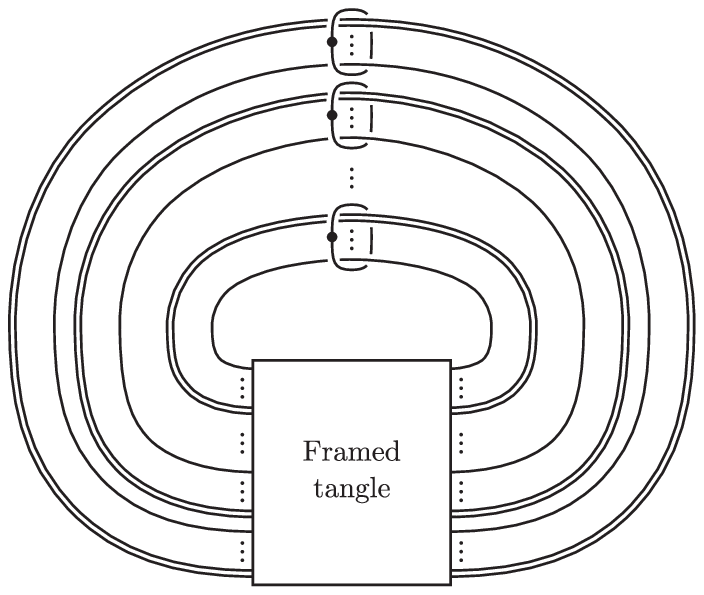}}
\end{Figure}

By the classical Alexander's argument, we can modify this diagram in order to make
the framed tangle inside the box into a framed braid. Moreover, by inserting a
certain number of kinks and enlarging them to form new braid strings, we can assume
that all the framings coincide with the blackboard ones. Figure \ref{diag2/fig}
shows the resulting diagram cut open in the upper part, after the 1-handles have
been isotoped to the lower part. 

\begin{Figure}[htb]{diag2/fig}{}{}\vskip-1mm
\centerline{\fig{Braided Kirby diagram}{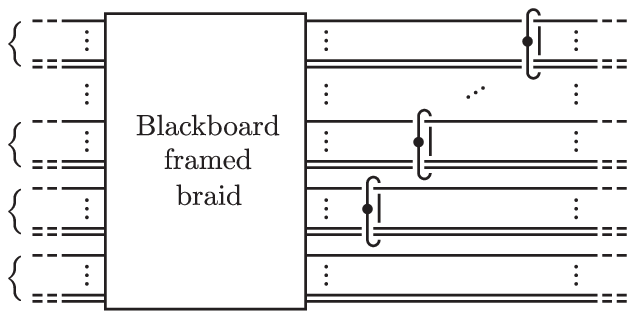}}
\end{Figure}

The rest of this section is devoted to show how the handle presentation of Figure
\ref{diag2/fig} can be converted into a simple covering $M \to B^4$ branched over an
orientable ribbon surface. We need first to specify some more details of
such handle presentation. Let $m$ and $n$ be respectively the number of 1-handles and
2-handles. We denote by $K^1_1, \dots, K^1_m$ the vertical trivial loops
representing the 1-handles in the diagram, indexed from left to right, and by
$K^2_1, \dots, K^2_n$ the braid components forming the attaching loops of the
2-handles, indexed according to their lowermost occurrence on the left, from bottom
to top. We assume the $K^2_i$'s counter-clockwise oriented. For any $i= 1, \dots,
n$, we call $s_i$ the number of strings of $K^2_i$ and we put $t_i = s_1 +
\dots + s_i$. As a notational convenience, we also put $t_0 = 0$. Moreover, $H^i_j$
will indicate the $i$-handle corresponding to $K^i_j$.

To begin with, we consider the simple branched covering of $B^4$ with $t_n + m + 1$
sheets numbered from $0$ to $t_n + m$, whose branching surface consists of the
trivial family of disjoint disks $D_1, \dots, D_{t_n + 2m} \subset B^4$ and whose
monodromy is given as follows: the disks $D_{t_{i-1}+1}, D_{t_{i-1}+2},
\dots,D_{t_i}$, that will be used for the 2-handle $H^2_i$, have respective
monodromies $(0\ t_{i-1}{+}1), (t_{i-1}{+}1 \ t_{i-1}{+}2), \dots, (t_i{-}1\ t_i)$;
the disks $D_{t_n + 2j - 1}$ and $D_{t_n + 2j}$, corresponding to the 1-handle
$H^1_j$, have the same monodromy $(0\ t_n{+}j)$.

\begin{Figure}[htb]{cov1/fig}{}{}
\centerline{\fig{Starting branching set}{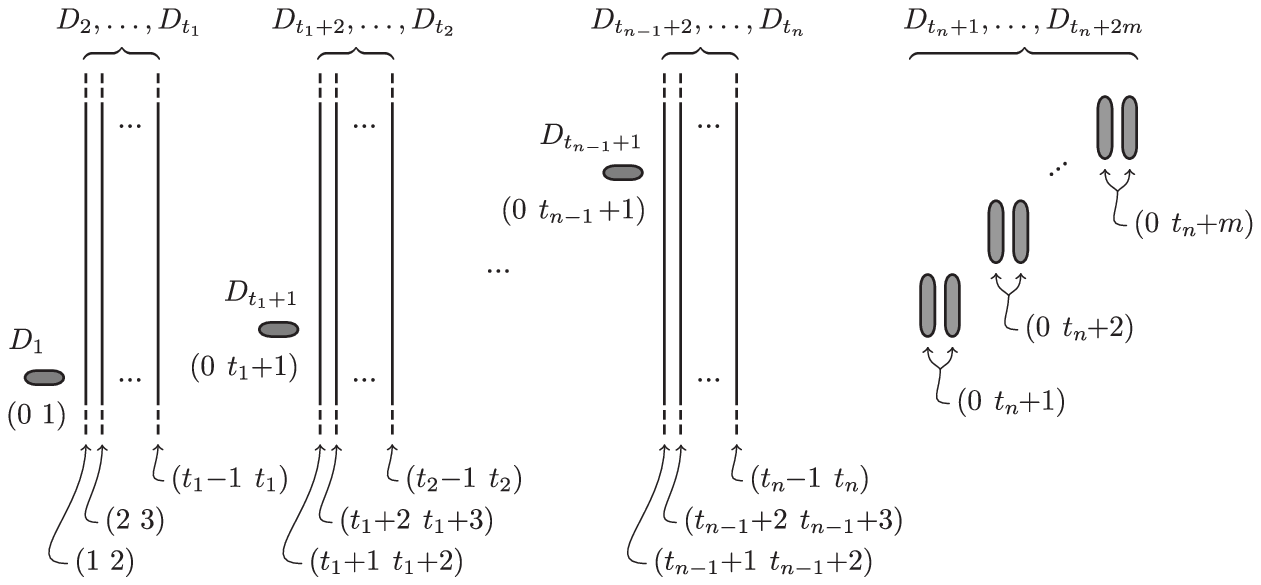}}
\end{Figure}

A diagram of these branching disks with their monodromies is shown in Figure
\ref{cov1/fig}. Here, the vertical lines stand for flat vertical disks, transversal
to the closed braid of Figure \ref{diag2/fig} in the upper part, where we cut open
it, so that each string meets all them once, from right to left in the order. There
is no such vertical disk for the 2-handles $H^2_i$ such that $K^2_i$ consists of
only one string, that this $s_i = 1$ and $t_{i-1}+1 = t_i$. Moreover, the disks
representing $D_{t_n + 2j - 1}$ and $D_{t_n + 2j}$ in the diagram are
$\epsilon$-displacements of the flat disk spanned by $K^1_j$ in Figure
\ref{diag2/fig}, hence the $K^2_i$'s cross these three disks in the same way, for
each $j = 1, \dots, m$. On the other hand, $D_{t_{i-1}+1}$ is a 2-disk expansion of a
small horizontal arc $C_i \subset K^2_i$ placed at the beginning (on the left in
Figure \ref{diag2/fig}) of the lowermost string of $K^2_i$, for each $i = 1,\dots,n$.

The covering manifold $M_1$ can be thought as $B^4 \cup H^1_1 \cup \dots \cup
H^1_m$. In fact, the disks $D_{t_n + 2j - 1}$ and $D_{t_n + 2j}$ give raise to the
1-handle $H^1_j$ formed by the sheet $t_n+j$, for each $j = 1, \dots, m$. All the
other branching disks induce stabilization by addition of trivial sheets. An outline
of $M_1$ (seen from the top) is drawn in Figure \ref{cov2/fig}. We identify $M_1$
with $B^4 \cup H^1_1 \cup \dots \cup H^1_m$ in such a way that the blackboard
framings relative to the Figures \ref{diag2/fig} and \ref{cov2/fig} coincide.

\begin{Figure}[htb]{cov2/fig}{}{}
\centerline{\fig{Starting branched covering}{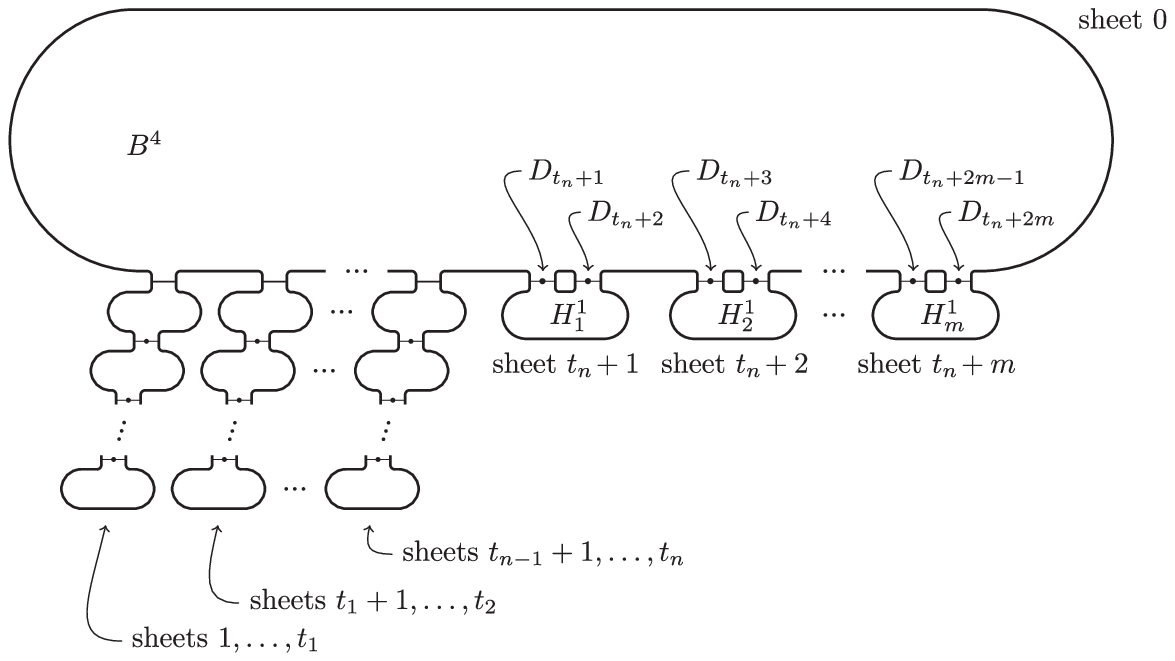}}
\end{Figure}

\medskip

Now, we modify the above branched covering $M_1 \to B^4$ to get the wanted
simple branched covering $M \cong M_1 \cup H^2_1 \cup \dots \cup H^2_n \to B^4$.

\break

Following Montesinos \cite{HM80} (see also \cite{IP02}), we realize the addition of
the 2-handles to $M_1$, by attaching an appropriate band $B_i$ to the branching
disk $D_{t_{i-1}+1}$, for each $i=1,\dots,m$. Namely, we define $B_i$ as a ribbon
band representing the blackboard framing along the arc $A_i = \Cl(K^2_i - C_i)$ in
Figure \ref{diag2/fig}.

This choice for the $B_i$'s works, since the following three properties are
satisfied (cf. \cite{HM80} or \cite{IP02}): 1) $A_i$ meets the branching disks only
at its endpoints, that belong to $\Bd D_{t_{i-1}+1}$; 2) the counterimage of $A_i$,
with respect to the covering, is the disjoint union of $t_n + m - 1$ arcs and a
simple loop $L_i = A'_i \cup A''_i \subset \Bd M_1$, where $A'_i$ and $A''_i$ are the
liftings of $A_i$ respectively starting in the sheets $0$ and $t_{i-1}+1$; 3) the
link $L_1 \cup \dots \cup L_n$ with the framings given by lifting the $B_i$'s is
equivalent (in $\Bd M_1$) to the link $K^2_1 \cup \dots \cup K^2_n$ with the
blackboard framings of Figure \ref{diag2/fig}.

Actually, property 1 holds by construction, while property 2 can be easily verified
by inspection, after observing that the product of the monodromies associated to the
vertical lines of Figure \ref{cov1/fig} taken from right to left is the permutation
$\pi = (1\ 2\ \dots\ t_1)(t_1{+}1\ t_1{+}2\ \dots\ t_2)\dots (t_{n-1}{+}1\
t_{n-1}{+}2\ \dots\ t_n)$.

\begin{Figure}[htb]{knot1/fig}{}{}\vskip0.5mm
\centerline{\fig{The arc $A_i$}{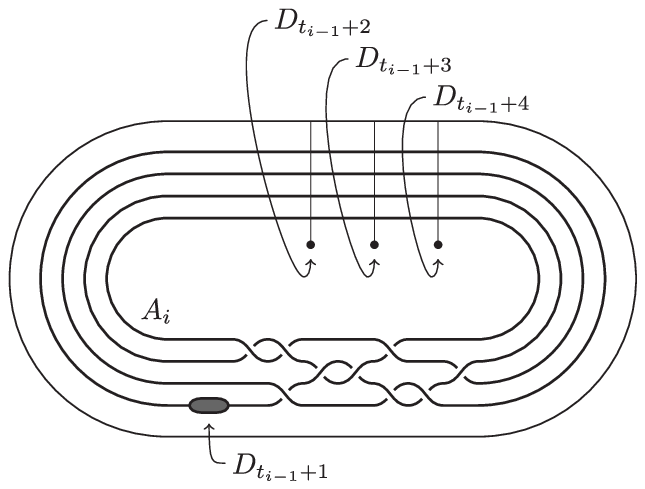}}
\end{Figure}

So, we are left with proving property 3. For the moment, we focus on a single braid
component $K^2_i$ disregarding the 1-handles. Figures \ref{knot1/fig} and
\ref{knot2/fig} describe respectively the arc $A_i$ and the loop $L_i$ for a braid
component $K^2_i$ with four strings, omitting the non-relevant branching disks and
sheets.

\begin{Figure}[htb]{knot2/fig}{}{}
\centerline{\fig{The loop $\widetilde K_i$}{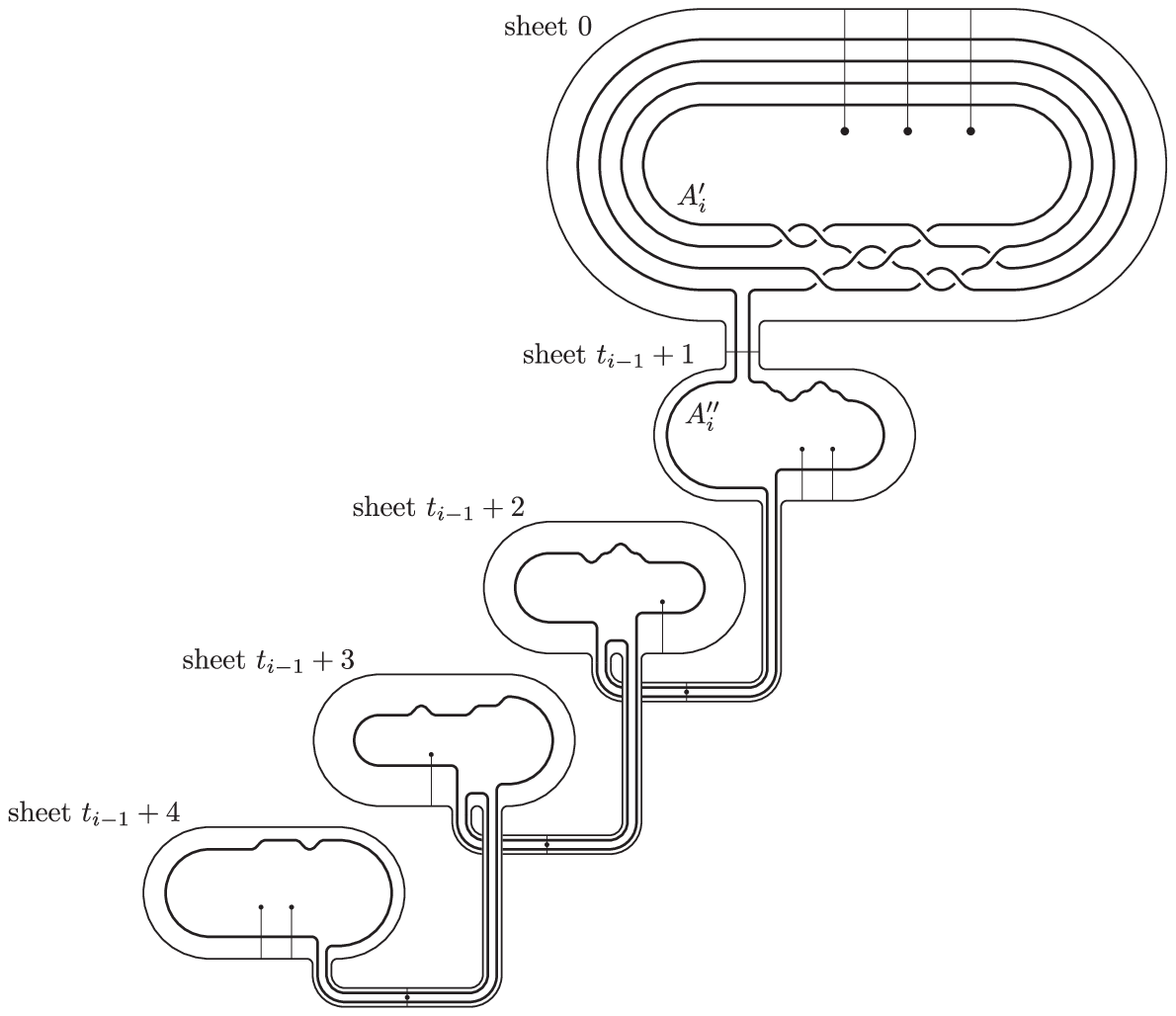}}
\end{Figure}

By considering again the permutation $\pi$, we immediately see that $A'_i$ is a copy
of $A_i$ entirely contained in the sheet $0$ of the covering, while $A''_i$ is a
trivial arc in the sheets $t_{i-1}+1, \dots, t_i$, consisting of one string in each
sheet (cf. Figure \ref{knot2/fig}). Hence, $L_i$ is isotopic to $K^2_i$.

Concerning the framing, we have that the blackboard framing along $A_i$ in Figure
\ref{cov1/fig} lifts to the blackboard framing along $L_i$ in Figure \ref{cov2/fig}
(cf. Figures \ref{knot1/fig} and \ref{knot2/fig}), which in turn is equivalent to the
blackboard framing of $K^2_i$ in Figure \ref{diag2/fig}. The last equivalence is due
to the fact that the isotopy between $L_i$ and $K^2_i$ can be assumed to be regular
with respect to the projection of figure \ref{cov2/fig}.

\medskip

\begin{Figure}[htb]{handle/fig}{}{}
\centerline{\fig{The 1-handles}{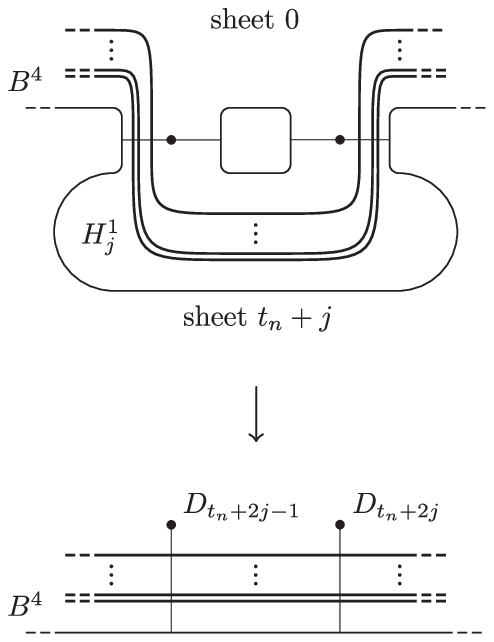}}
\end{Figure}

At this point, we observe that all the $K^2_i$'s can be considered simultaneously,
since the $L_i$'s interacts only in the sheet $0$. Hence $L_1 \cup \dots \cup L_n$
and $K^2_1 \cup \dots \cup K^2_n$ are isotopic as blackboard framed links.

\medskip

Finally, let us take into account the 1-handles. Figure \ref{handle/fig} shows how
the crossings of the $A_i$'s with the projections of $\Int D_{t_n + 2j -1}$ and $\Int
D_{t_n + 2j}$ lifts to passages of the $A'_i$'s through the 1-handle $H^1_j$. In
particular, we have that no extra twist is added neither in the link nor in the
framings. Then, the presence of the 1-handles does not affect our reasoning in any
way, except for the fact that the arcs $A'_i$ are no longer contained only in the
sheet $0$, but they traverse also the sheets $t_n + 1,\dots, t_n + m$ forming the
1-handles.

\break

\begin{Figure}[htb]{cov3/fig}{}{}\vskip2mm
\centerline{\fig{Special branched covering}{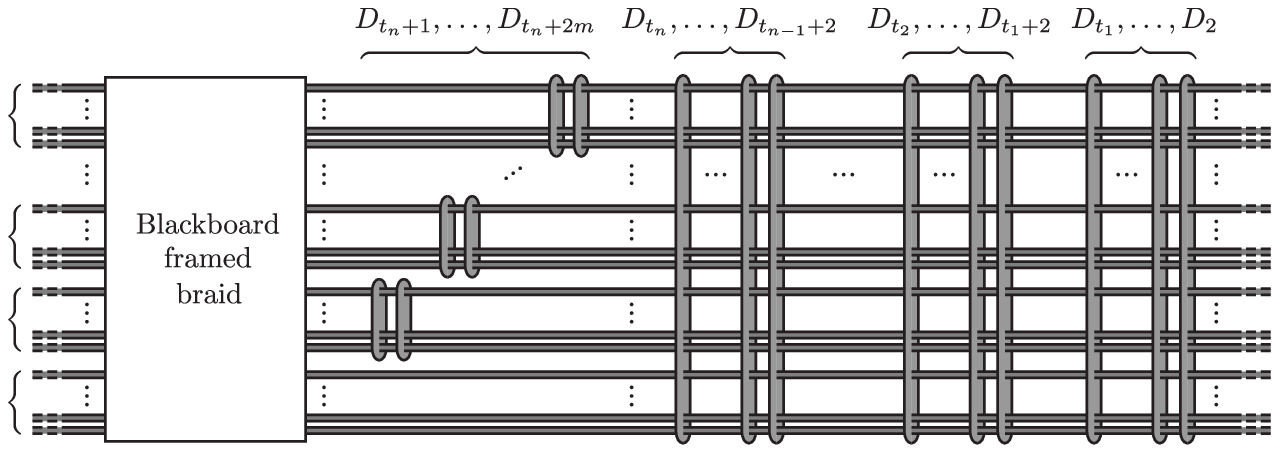}}
\end{Figure}

A diagram of the resulting branching surface (cut open as above) is outlined in
Figure \ref{cov3/fig}. Here, the framed braid is the same of Figure \ref{diag2/fig},
while the vertical disks are the ones of Figure \ref{cov1/fig}, apart from different
order, due to the sliding of the $D_k$'s with $k \leq t_n$ from the upper part of the
diagram (cf. Figure \ref{knot1/fig}) to the lower one.


\section{Getting the universal surface\label{surface/sec}}

We begin this section, by explaining how the covering moves given in Section
\ref{moves/sec} can be used to symmetrize the branching surface of Figure
\ref{cov3/fig}.

Firstly, we modify any positive (resp. negative) crossing along the braid inside the
box as described in the top (resp. bottom) part of Figure \ref{move5/fig}. In both
cases, we perform eight Moves 3 (cf. Figure \ref{move3/fig}) and then we isotope
some of the resulting vertical disks. Then, we make all such crossings into ribbon
intersections, by stabilization (followed by suitable isotopy) and crossing change,
as shown in Figure \ref{move6/fig} (of course, the covering degree $d$ must be
dynamically updated after each stabilization). We leave to the reader to check that
the monodromies of the two bands forming each crossing are really distinct but not
disjoint, like in Figure \ref{move6/fig}.

\begin{Figure}[htb]{move5/fig}{}{}
\centerline{\fig{Move 5}{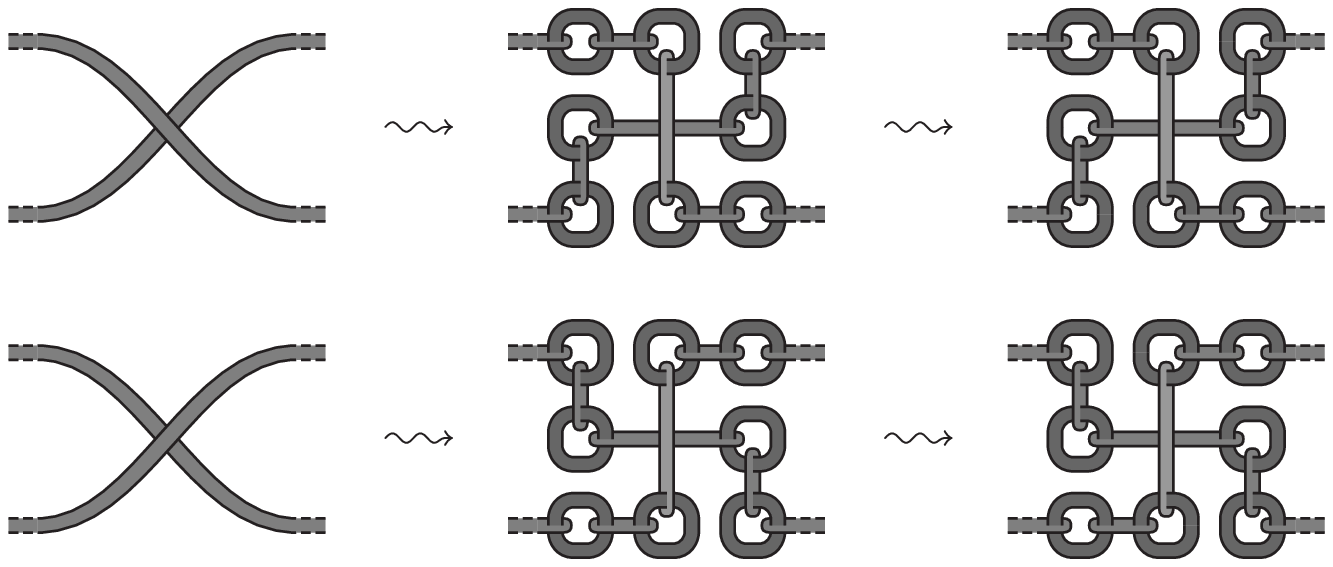}}
\end{Figure}

\begin{Figure}[htb]{move6/fig}{}{}\vskip1mm
\centerline{\fig{Move 6}{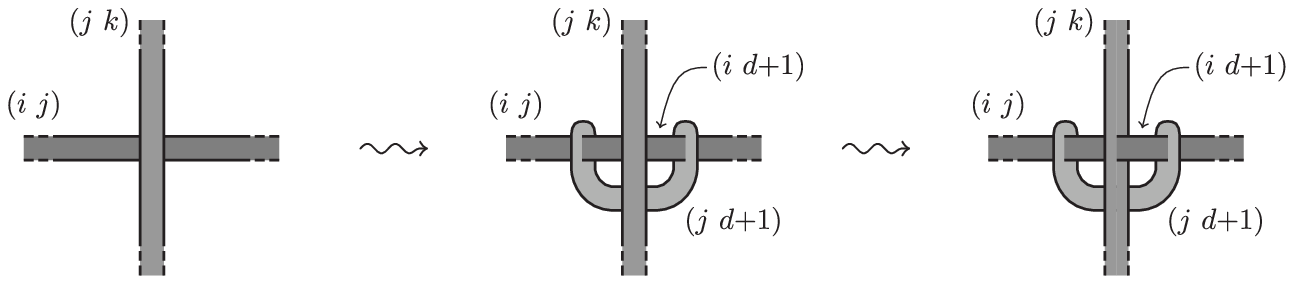}}
\end{Figure}

Secondly, we apply our Move 4 (cf. Figure \ref{move4/fig}) to the ribbon
intersections we have just obtained, apart from the ones formed by the
stabilizing disks. Moreover, we do the same on all the ribbon intersections which
appear in Figure \ref{cov3/fig} outside the box. Also in this case, we leave to the
reader to verify that the involved monodromies are distinct.

At this point, our diagram consists of: small annuli centered at some vertices of a
rectangular grid; a certain number of horizontal and vertical bands running along
some edges of the same grid; small stabilizing disks as in Figure \ref{move6/fig}
around some of the annuli. We emphasize that the bands do not form any ribbon
intersection or crossing with each other. 

\begin{Figure}[b]{surf1/fig}{}{}
\centerline{\fig{Surface 1}{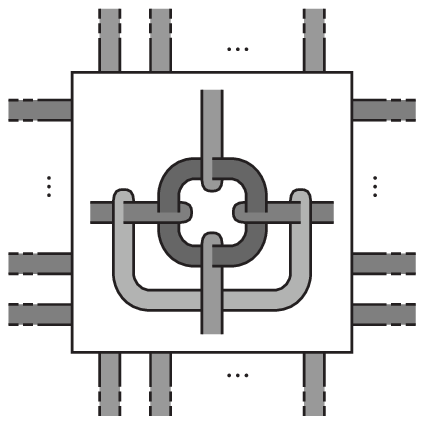}}
\end{Figure}

Such a diagram can be easily completed to get the one depicted in Figure
\ref{surf1/fig}, where top and bottom ends are assumed to be trivially joined by
bands passing in front of the ones already connecting left and right ends, and the
pattern decorating the box has to be replicated at each potential crossing between
the entering horizontal and vertical bands. Namely, it suffices to break the bands
which take more that one grid edge, by using Move 3, and then to insert fake
branching components (labelled with the identity) in the lacking places. Of course,
top and bottom grid lines have to be considered as if they were adjacent.

\medskip

Now, thinking of $B^4$ as $B^2 \times B^2 \subset \Bbb C^2$, we place the diagram of
Figure \ref{surf1/fig} in the torus $S^1 \times S^1$, in such a way that the
rotations $r_1: (z_1,z_2) \mapsto (e^{2\pi i/n_1}z_1,z_2)$ and $r_2: (z_1,z_2)
\mapsto (z_1,e^{2\pi i/n_2}z_2)$ permute respectively the rows and the columns of
the $n_1 \times n_2$ pattern matrix inside the box (cf. \cite{T82} and
\cite{HLMW87}).

\medskip

Then, we compose the branched covering represented by the diagram with the quotient
by the action of $r_2$, to get a new branched covering $M \to B^4$, whose branching
surface is given in Figure \ref{surf2/fig}. Here, the rightmost disk is the
branching surface of the quotient, while the box contains a $n_1 \times 1$ pattern
matrix, which is the quotient of the one of Figure \ref{surf1/fig}. 

\begin{Figure}[htb]{surf2/fig}{}{}
\centerline{\fig{Surface 2}{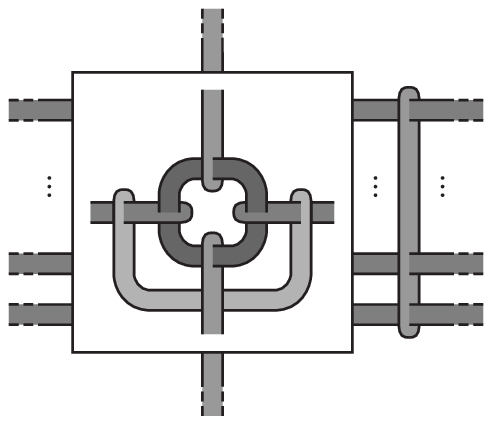}}
\end{Figure}

Breaking the rightmost disk into $n_1$ disks, by Move 3 once again, and adding
another fake branching annulus between top and bottom, we get the diagram of Figure
\ref{surf3/fig}, which can be still assumed $r_1$-invariant.

\begin{Figure}[htb]{surf3/fig}{}{}\vskip1mm
\centerline{\fig{Surface 3}{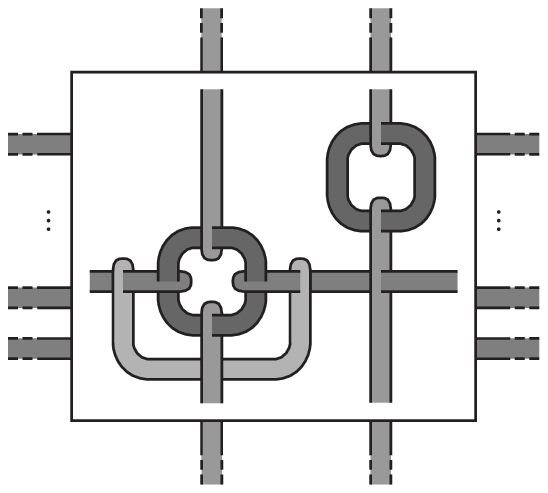}}
\end{Figure}

Finally, we quotient by the action of $r_1$, in order to get the diagram of Figure
\ref{surf4/fig}, where the branching disk of this last quotient is the horizontal
one. It is worth remarking that, by quotienting directly the diagram of Figure
\ref{surf2/fig}, one would get a singular point in the surface, due to the
transversal intersection between the branching disks of the two quotients. This is
the motivation for passing through the diagram of Figure \ref{surf3/fig}.

\break

\begin{Figure}[htb]{surf4/fig}{}{}
\centerline{\fig{Surface 4}{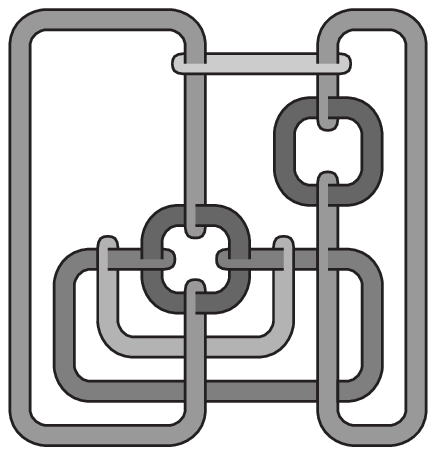}}
\end{Figure}

Clearly, Figure \ref{surf4/fig} already represents a universal orientable ribbon
surface. However, we conclude this section by simplifying a little bit such
universal surface. The intermediate steps of this simplification are described in
the following Figure \ref{surf5/fig}. 

\begin{Figure}[htb]{surf5/fig}{}{}\vskip0.5mm
\centerline{\fig{Surface 5}{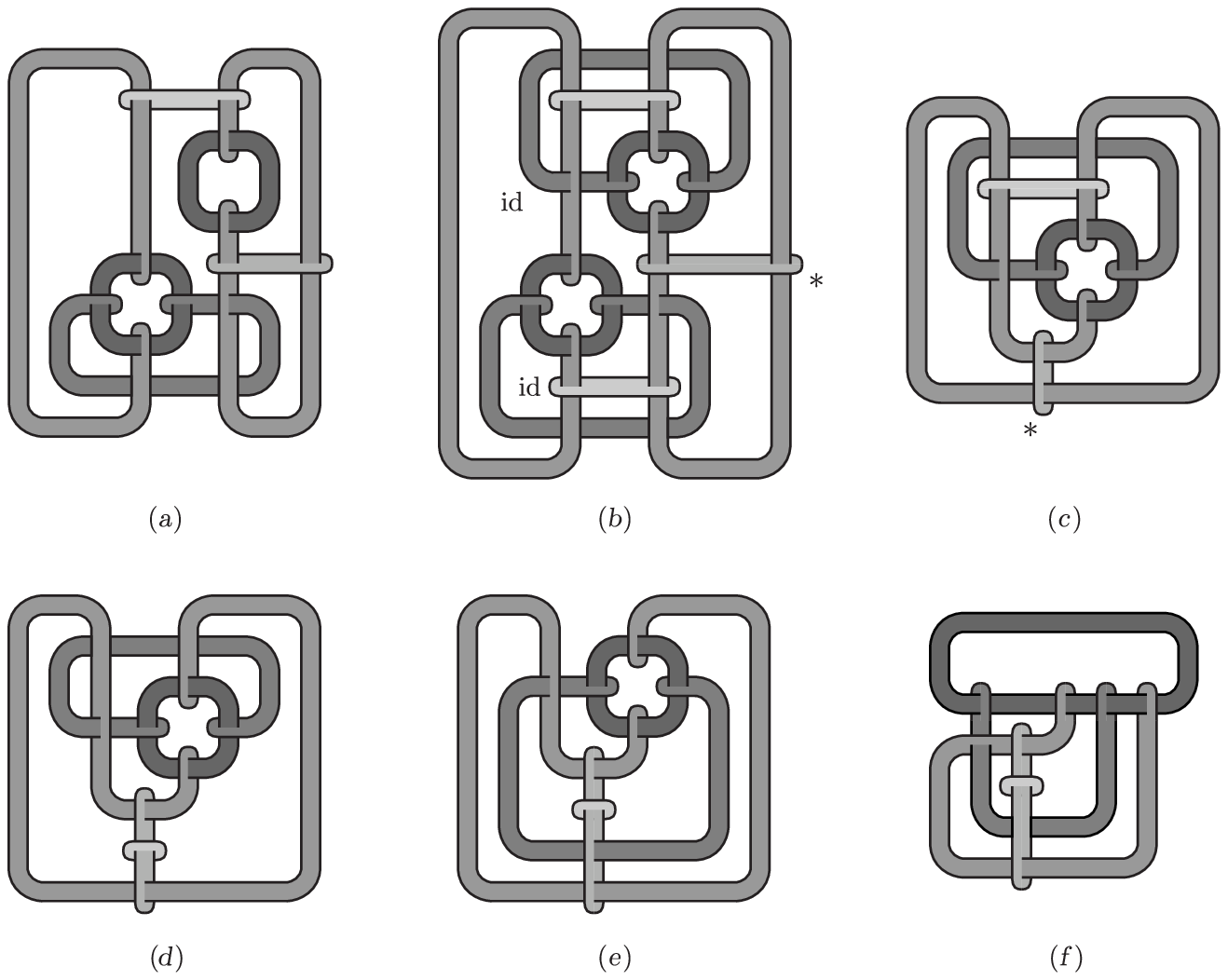}}
\end{Figure}

We start with the surface in \(a), which is isotopically equivalent to the one
of Figure \ref{surf4/fig}. Then, we add the two fake branching components labelled
with the identity in \(b), in order to make the surface symmetric with respect
the center of the diagram (the disk marked with the asterisk can be thought as the
fixed point set of the symmetry). In \(c) we see the branching surface of the
composition of our branched covering with that symmetry. Of course, this surface is
still universal, but it has two components less than before. The surfaces in \(d),
\(e) and \(f) are all obtained by isotopy.

\break

The simplified universal surface is depicted in Figure \ref{univ/fig}. To get it, we
isotoped once again the last surface of Figure \ref{surf5/fig} just for pictorial
reasons.

\begin{Figure}[htb]{univ/fig}{}{}\vskip1mm
\centerline{\fig{Universal surface}{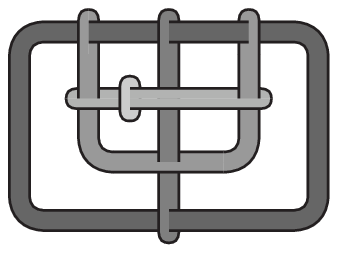}}
\end{Figure}

\section{Concluding remarks and questions\label{remarks/sec}}

First of all, we observe that the universal surface of Figure \ref{univ/fig}
consists of one annulus and four disks, all trivially embedded in $B^4$. Moreover,
disregarding the annulus, the four disks can be separated and isotoped in a
symmetric position, so that they are cyclically permuted by a rotation of $\pi/2$
radians. Then, the quotient by the action of the rotation gives us a new universal
ribbon surface with only three components, one annulus and two trivial disks.
Unfortunately, the isotopy and the quotient force the annulus to wrap around the
disks in a very unpleasant way and this makes resulting surface likely useless.
Nevertheless, we know that the number of components can be reduced to three.
However, the following question makes sense.

\begin{question}
Is there a ``reasonable'' universal (possibly non-orientable) surface in
$B^4$ with less that five components?
\end{question}

Even more, there is no reason to believe that three is the minimum number of
components of a universal surface in $B^4$. In fact, it can be easily proved, by
using signature, that there is no connected universal surface in $S^4$ (cf.
\cite{V84} and \cite{IP02}), but the same argument does not work in $B^4$. 
So, here is our second question.

\begin{question}
Does there exist any connected universal surface in $B^4$?%
\end{question}

\medskip

On the other hand, at the cost of some more components, one could modify the
construction carried out in Section \ref{surface/sec}, in order to get a different
universal surface, such that branched covering with all the branching indices equal
to $2$ would suffice for our representation theorem. The only branching indices
bigger than $2$ coming into that construction are indeed due to the rotations $r_1$
and $r_2$. Namely, the branching indices over the two disks fixed by such rotations
(cf. Figures \ref{surf2/fig} and \ref{surf4/fig}) are respectively $n_1$ and $n_2$.
By the merging principle, each of these disks can be replaced by a pair of parallel
disks labelled with suitable products of disjoint transpositions (the same argument
used in \cite{HLMW87} for the 3-dimensional case applies here). In this way, all the
branching indices are reduced to $2$.

\medskip

Figure \ref{surf6/fig} \(b) shows a braided version of our universal surface.
It has been obtained by applying the Rudolph's braiding algorithm (cf. \cite{R83})
to the surface in the part \(a) of the same figure, which is isotopic to the one of
Figure \ref{univ/fig}. 

\begin{Figure}[htb]{surf6/fig}{}{}
\centerline{\fig{Surface 6}{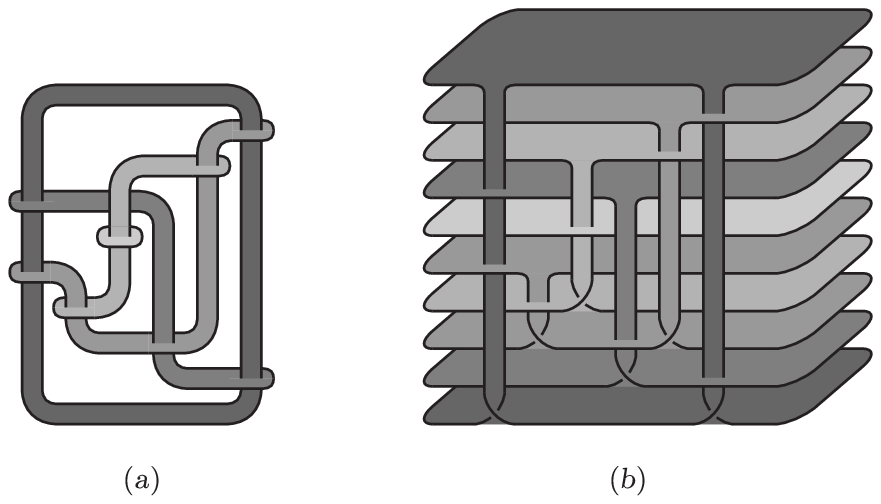}}
\end{Figure}

Such way of seeing a universal surface is quite interesting, due to that fact that
any covering $M \to B^4$ branched over a braided surface naturally induces a
topological Lefschetz fibration $M \to B^2$ (cf. \cite{LP01}). 
In particular, if the braided surface is positive (quasi-positive in the Rudolph's
terminology), that is all the bands between the sheets are positively twisted, then
also the induced fibration is positive and $M$ is a bounded Stein surface.
One of the main ingredients in the proof of this fact is the Rudolph's theorem that
positive braided surfaces are complex analytic (cf. \cite{R83a}).

Since not all the 4-manifold considered here are Stein, no universal ribbon surface
in $B^4$ can be isotopically equivalent to a positive braided surface. In Figure
\ref{surf6/fig} \(b) we see that only three of the six bands are positively twisted.

However, it has been proved in \cite{LP01} that any compact Stein surface with
boundary is a covering of $B^4$ branched over a positive braided surface. Then, the
following question naturally arises. 

\begin{question}
Does there exists a positive braided surface in $B^4$ which is universal for compact
Stein surfaces with boundary?
\end{question}

\medskip

Finally, some trivial remarks about universal links. Obviously, the boundary of any
universal surface in $B^4$ is a universal link in $S^3$. But, it is likely false that
any universal link in $S^3$ is the boundary of a universal surface in $B^4$.
Actually, it is not clear at all whether any universal link bounds some surface in
$B^4$ allowing us to give a covering presentation of any closed orientable
3-manifold as the boundary of a 4-manifold. For example, we don't know what happens
in the simplest case of the Borromean rings. So, we conclude with the following
question.

\begin{question}
What universal links in $S^3$ bound a universal surface in $B^4$?
\end{question}

\break

\bibitem{H74}	H.M. Hilden, {\sl Every closed orientable 3-manifold is a 3-fold
branched covering space of $S^3$}, Bull. Amer. Math. Soc. {\bf 80} (1974),
1243--1244.

\bibitem{HLM83}	H.M. Hilden, M.T. Lozano and J.M. Montesinos, {\sl The Whitehead
link, the Borromean rings and the knot 946 are universal}, Collect. Math. {\bf 34}
(1983), 19--28.

\bibitem{HLM85}	H.M. Hilden, M.T. Lozano and J.M. Montesinos, {\sl Universal knots},
Lecture Notes in Math. {\bf 1144}, Springer-Verlag 1985, 25--59.

\bibitem{HLM85'}	H.M. Hilden, M.T. Lozano and J.M. Montesinos, {\sl On knots that are
universal}, Topology {\bf 24} (1985), 499--504.

\bibitem{HLMW87}	H.M. Hilden, M.T. Lozano, J.M. Montesinos and W.C. Whitten, {\sl On
universal groups and three-manifolds}, Invent. Math. {\bf 87} (1987), 441--456.

\bibitem{HM80} H.M. Hilden and J.M. Montesinos, {\sl Lifting surgeries to branched
covering\break spaces}, Trans. Amer. Math. Soc. {\bf 259} (1980), 157--165.

\bibitem{Hi74}	U. Hirsch, {\sl \"Uber offene Abbildungen auf die 3-Sph\"are}, Math.
Z. {\bf 140} (1974), 203--230.

\bibitem{IP02} M. Iori and R. Piergallini, {\sl 4-manifolds as covers of $S^4$
branched over non-singular surfaces}, preprint.

\bibitem{K89} R. Kirby, {\sl The topology of 4-manifolds}, Lecture Notes in
Mathematics {\bf 1374}, Springer-Verlag 1989.

\bibitem{LP01} A. Loi and R. Piergallini, {\sl Compact Stein surfaces with boundary
as branched covers of $S^4$}, Invent. math. {\bf 143} (2001), 325--348.

\bibitem{M74}	J.M. Montesinos, {\sl A representation of closed, orientable
3-manifolds as 3-fold branched coverings of $S^3$}, Bull. Amer. Math. Soc. {\bf 80}
(1974), 845--846.

\bibitem{M78} J.M. Montesinos, {\sl 4-manifolds, 3-fold covering spaces and
ribbons}, Trans. Amer. Math. Soc. {\bf 245} (1978), 453--467.

\bibitem{M83}	J.M. Montesinos, {\sl Representing 3-manifolds by a universal
branching set}, Proc. Camb. Phil. Soc. {\bf 94} (1983), 109--123.

\bibitem{M85}	J.M. Montesinos, {\sl A note on moves and irregular coverings of
$S^4$}, Contemp. Math. {\bf 44} (1985), 345--349.

\bibitem{MP01} M. Mulazzani and R. Piergallini, {\sl Lifting braids}, Rend. Ist.
Mat. Univ. Trieste {\bf XXXII} (2001), Suppl. 1, 193--219.

\bibitem{P95} R. Piergallini, {\sl Four-manifolds as $4$-fold branched covers of
$\S^4$}, Topology {\bf 34} (1995), 497-508.

\bibitem{R83} L. Rudolph, {\sl Braided surfaces and Seifert ribbons for closed
braids}, Comment. Math. Helvetici {\bf 58} (1983), 1--37.

\bibitem{R83a} L. Rudolph, {\sl Algebraic functions and closed braids}, Topology {\bf
22} (1983),\break 191--202.

\bibitem{T82} W. Thurston, {\sl Universal links}, preprint 1982.

\bibitem{V84} O.Ja. Viro, {\sl Signature of branched covering},
Trans. Mat. Zametki {\bf 36} (1984), 549--557.


\thebibliography{[00]}

\bibitem{Ap03} N. Apostolakis, {\sl On 4-fold covering moves}, Algebraic \&
Geometric Topology {\bf 3} (2003), 117--145. 

\bibitem{BE79}	I. Bernstein and A.L. Edmonds, {\sl On the construction of branched
coverings of low-dimensional manifolds}, Trans. Amer. Math. Soc. {\bf 247} (1979),
87--124.

\bibitem{BP03} I. Bobtcheva and R. Piergallini, {\sl On $n$-fold covering moves},
preprint.


\bibitem{Hr03} F. Harou, {\sl Description en terme de rev\^etements simples de
rev\^etements ramifi\'es de la sph\`ere}, preprint.
\end{document}